\documentstyle{article}
\newtheorem{tm}{Theorem}[subsection]
\newtheorem{lm}[tm]{Lemma}
\newtheorem{pr}[tm]{Proposition}
\newtheorem{rmk}[tm]{Remark}
\newtheorem{cor}[tm]{Corollary}
\newtheorem{ex}[tm]{Example}

\newtheorem{??}[tm]{Question}

\font\tenmsb=msbm10
\font\sevenmsb=msbm7
\font\fivemsb=msbm5

\newfam\msbfam
\textfont\msbfam=\tenmsb
\scriptfont\msbfam=\sevenmsb
\scriptscriptfont\msbfam=\fivemsb
\def\Bbb#1{{\fam\msbfam #1}}

\font\teneufm=eufm10
\font\seveneufm=eufm7
\font\fiveeufm=eufm5
\newfam\eufmfam
\textfont\eufmfam=\teneufm
\scriptfont\eufmfam=\seveneufm
\scriptscriptfont\eufmfam=\fiveeufm
\def\frak#1{{\fam\eufmfam\relax#1}}

\newcommand\ci{\cite}
\newcommand\s{\sigma}
\newcommand\rat{{\Bbb Q}}
\newcommand\comp{{\Bbb C}}
\newcommand\real{{\Bbb R}}
\newcommand\zed{{\Bbb Z}}
\newcommand\nat{{\Bbb N}}
\newcommand\pn[1]{{\Bbb P}^{#1}}

\newcommand\blacksquare{{\hspace*{\fill} $\Box$}} 
\newcommand\odix[1]{ {\cal O}_{#1} }

\newcommand\e{\epsilon}

\newcommand{\SS}{\frak S}

\newcommand\D{\Delta}

\newcommand{\bnu}{{\tilde{\nu}}}
\newcommand{\ba}[1]{{\bf {#1}}}

\newcommand{\sbcpt}{ _z\D^{(\bnu)}_{\beta} (\epsilon) }

\newcommand{\XA}{X^{(\ba{a})}}

\newcommand{\Xa}{X^{(n)}_{(\ba{a})}}

\newcommand{\n}{\noindent}
 
\newcommand{\be}{\begin{equation}}
\newcommand{\ee}{\end{equation}}
\newcommand{\ea}{\end{array}}
\def\lorw{\longrightarrow}

\def \H#1{ X^{[#1]}}
\def\ee{\`{e}}

\renewcommand{\ss}[1]{\medbreak{\bf#1.}\kern4pt}

\begin{document}

\centerline{\bf{\huge The Douady space of a complex surface}}

%\maketitle

\vspace{0.7cm}

\centerline{\Large
Mark Andrea A.  de Cataldo\footnote[1]{
Partially supported by N.S.F. Grant DMS 9701779 and 
by an A.M.S. Centennial Fellowship.} 
and  Luca Migliorini\footnote[2]{Supported
by M.U.R.S.T. funds and by SFB 237.}}

\bigskip
\centerline{November 24, 1998}

\medskip
%\centerline{\bf {\huge Please do not distribute}}

%\date{}

\vspace{1cm}

\begin{abstract}
We prove that a standard realization 
of the direct image complex via the so-called
Douady-Barlet morphism associated with a smooth complex analytic surface
 admits a natural decomposition in the form of an injective
quasi-isomorphism of complexes.
This is a more precise form of a special case of the decomposition theorems
of Beilinson-Bernstein-Deligne-Gabber and M. Saito.

The proof hinges on the special case of the bi-disk in the complex affine plane
where  we   make explicit use of  a construction 
of
Nakajima's 
and of the corresponding representation-theoretic interpretation
foreseen by Vafa-Witten. 
%to 
 %derive a simple proof of a special case of the  the so-called G\"ottsche %Formula, i.e. a formula for the Betti numbers of the Douady spaces of a bi-disk %,and  obtain canonical and geometrically meaningful bases for the rational %cohomology of these spaces. These
%bases are then linked, in the  case of a complex analytic surface, 
%to the natural stratification of the Douady-Barlet morphisms. 

%Our approach underlines the importance of the representation
%theoretic side, gives new and elementary proofs of G\"ottsche formula and of %the decomposition theorem in the algebraic case, and extends the theory
%to the complex analytic case.

Some consequences of the decomposition theorem:
G\"ottsche Formula holds for complex surfaces; interpretation
 of the rational  cohomologies
of  Douady spaces as a kind of Fock space; new proofs of results of Brian\c{c}on and Ellingsrud-Stromme
on punctual Hilbert schemes; computation of the mixed Hodge structure 
of the Douady spaces 
in the K\"ahler case. We also derive a natural  connection with Equivariant K-Theory for which,    in the case of algebraic surfaces,
Bezrukavnikov-Ginzburg have proposed a different approach.

\end{abstract}

\tableofcontents

\section{Introduction}
Let $X$ be a smooth complex analytic surface, $n$ be a non-negative integer
and $X^{[n]}$ be the Douady space of zero-dimensional analytic subspaces of $X$ of length $n$. The spaces $X^{[n]}$ are $2n$-dimensional complex manifolds.
If $X$ is algebraic, then they are the usual Hilbert schemes.

\noindent
These  manifolds have been intensively studied from a ``local"
perspective; see \ci{fo}, \ci{ia}, \ci{ia2}, \ci{ia3},
\ci{br},  \ci{e-s}, \ci{e-s2}; this list is by no means complete.

\medskip
Let $X$ be algebraic. G\"ottsche Formula (see Theorem \ref{go}.(\ref{gofock}))
gives the Betti numbers of $X^{[n]}$ in terms of the ones of $X$.

\noindent
Vafa and Witten \ci{va-wi} remarked that G\"ottsche Formula is the character
formula of the standard {\em irreducible} representation
 of a certain infinite dimensional
super Lie algebra ${\cal H}(X)$, called the Heisenberg/Clifford algebra,  modeled on the rational cohomology
$H^*(X)$ of $X$. As a consequence, the vector space ${\Bbb H}(X):
= \oplus_{n\geq 0}{H^*( X^{[n]}  )}$ can be seen, abstractly, as an
{\em irreducible}  highest weight ${\cal H}(X)$-module.

\medskip
Motivated by this remark, Nakajima \ci{nak} realized {\em geometrically}, and 
for every complex surface $X$, the space ${\Bbb H}(X)$ as a 
highest weight ${\cal H}(X)$-module by introducing certain 
operators acting on ${\Bbb H}(X)$ induced by correspondences in the products of Douady spaces.  Grojnowski has announced similar results  \ci{gr}.

\noindent
In this context, G\"ottsche Formula  Theorem \ref{go}.(\ref{gofock}) 
becomes the statement
that  if $X$ is algebraic, then the geometrically realized action 
of  
${\cal H}(X)$ on ${\Bbb H}(X)$ is irreducible. In particular,
${\Bbb H}(X)$ can be re-built, at least in principle,  from $H^*(X)$
by means of explicit geometric operators.

\medskip
This new and beautiful structure of the Hilbert schemes $X^{[n]}$ ($X$ algebraic) 
emerges because all values of $n$ have been considered {\em simultaneously}.

\bigskip
All known  proofs of G\"ottsche Formula, for $X$ algebraic, rely,
to start with, on results of 
Brian\c{c}on \ci{br} and Ellingsrud and Stromme \ci{e-s} on punctual
Hilbert schemes. Each individual proof then relies on either
Deligne's solution to the Weil Conjectures \ci{go}, on Beilinson-Bernstein-Deligne-Gabber Decomposition Theorem
for perverse sheaves \ci{go-so}, or on Deligne's theory of
mixed Hodge structures 
\ci{ch}.

\bigskip
In this paper we develop a new approach by taking the ${\cal H}(X)$-action
on ${\Bbb H}(X)$ as a starting point.

\noindent
Our approach, inspired by Nakajima's \ci{nak}, but otherwise
self-contained, underlines the importance of the 
representation-theoretic side, gives an 
elementary proof of the relevant decomposition theorem and of G\"ottsche 
Formula, and it extends the theory to the complex analytic case.

\noindent
In  particular, {\em all} the results mentioned above, which hold 
 in the algebraic case, are proved here to hold, more generally,
 in the analytic context. They are all consequences 
of our analysis  of the special case $X=\Delta$ the bi-disk in $\comp^2$ and
of our Decomposition Theorem \ref{decel}. 

\noindent
Moreover, we establish a new, natural relation with
Equivariant $K$-Theory.

\smallskip

This analytic viewpoint would seem to be more natural.

\medskip
Nakajima \ci{nak} asks whether it is possible to
extend the picture drawn in the algebraic case to differentiable four-manifolds
and in this paper we show that this can be done for complex analytic surfaces.
 
\medskip
The analysis
 builds on the special case $X=\D$. Here the role of the usual Heisenberg algebra ${\cal H}(\D)$ is paramount: we determine an explicit, natural
and geometrically meaningful basis for ${\Bbb H}(\D)$. The case $X=\comp^2$
is analogous and we recover results in \ci{e-s} with the basic difference that
our basis differs from the one in  \ci{e-s1}. We exploit the geometric meaning of our  basis for ${\Bbb H}(\D)$ and
prove, {\em for every} complex analytic surface $X$, a precise form
of the Decomposition Theorem
for the Douady-Barlet morphism $\pi: X^{[n]} \to X^{(n)}$ in the form
of an injective quasi-isomorphism of complexes. The relevant morphism
in \ci{b-b-d}, used in \ci{go-so} to prove G\"ottsche Formula, is in a  derived  category and it is not uniquely determined.

\noindent
G\"ottsche Formula
and the irreducibility of ${\Bbb H}(X)$ follow formally, and so does
the determination, first obtained in \ci{go-so} using Saito's theory of mixed Hodge modules, of the 
 mixed Hodge structure of $X^{[n]}$ in the algebraic (or ``K\"ahler") case.
We also re-obtain results in \ci{br} and in  \ci{e-s} on punctual Hilbert schemes.
Another  easy consequence of our Decomposition
Theorem is the construction of 
a natural additive isomorphism between the rational $\SS_n$-Equivariant $K$-Theory
of $X^n$ (here $\SS_n$ acts on $X^n$ by permuting the factors)
and the rational $K$-Theory of $X^{[n]}\, :$
$K_{\SS_n} (X^n) \otimes \rat \simeq K(X^{[n]})\otimes \rat$.
Bezrukavnikov and Ginzburg \ci{b-g} have proposed to construct, 
in the algebraic case,
a   different natural map.
Our motivation and proof are different from theirs.
 
\bigskip 
The paper is organized as follows.
In Section \ref{dsbd}
we give a new, explicit, construction of $X^{[n]}$
and of the Douady-Barlet morphism $\pi: X^{[n]} \to X^{(n)}$.
The building blocks are the Douady spaces for the bi-disk $\D^{[m]}$, $m\leq n$.
These, in turn, are constructed using the ``toy model" of Nakajima \ci{naknotes}.
For $X=\D$, $\comp^2$, we use the ${\cal H}(X)$-action  on ${\Bbb H}(X)$ to compute the Betti numbers of
$X^{[n]}$ and  determine a canonical basis for 
${\Bbb H}(X)$, Theorem \ref{pobd}.  We use this result to re-prove
Ellingsrud-Stromme formula for the Betti numbers of punctual Hilbert schemes
of a surface, Corollary \ref{bettibriancon}, and the irreducibility result of
Brian\c{c}on, Corollary \ref{brianc}.

\noindent
Section 
\ref{douadybarlet} is preparatory 
for the Decomposition Theorem. We study the Douady-Barlet morphism and its natural stratification. The normalizations of the closures of the strata in
$X^{(n)}$ play a basic role. This section discusses the natural identifications
which we obtain  between objects on these normalizations, on $X^{(n)}$ and on
$X^{[n]}$. We then define a certain injective morphism
of complexes $\Psi$; see Proposition \ref{prpsi}.

\noindent
In Section \ref{decompo} we prove our Decomposition Theorem,
Theorem \ref{decel}: $\Psi$ is a quasi-isomorphism. 

\noindent
In Section \ref{conseq} we deduce formal consequences of Theorem
\ref{decel} which are new in the analytic case:
Corollary \ref{isr} ($R^q\pi_* \rat_{X^{[n]}}$), Theorem \ref{leray} (degeneration of the Leray spectral sequence), Theorem \ref{go}
(G\"ottsche Formula), 
Corollary \ref{mark} (Euler numbers of $X^{[n]}$), Theorem \ref{naka} (irreducibility
of ${\Bbb H}(X)$ as a highest weight ${\cal H}(X)$-module), Theorem \ref{mhsxn} (mixed Hodge structure
of $X^{[n]}$).
Finally we prove Theorem \ref{eqco} ($K_{\SS_n} (X^n) \otimes \rat \simeq K(X^{[n]})\otimes \rat$).

\bigskip
\noindent
{\bf Acknowledgments.}
This work has been partly inspired by H. Nakajima's papers \ci{nak} and \ci{naknotes}; we heartly thank him for useful correspondences.

\noindent
W. Wang's ideas and enthusiasm have influenced our work.

\noindent
Parts of this work have been carried out at the Max-Planck-Institut f\"ur Mathematik in Bonn, at R\"uhr-Universit\"at-Bochum and at Harvard University.

\noindent
We want to thank all the participants in the seminars on Nakajima's work which have taken place at the M.P.I. and at the University of Rome ``La Sapienza,"
especially,  M. Furuta, M. Kim and
A. Raina. We also want to thank  T. Iarrobino, A. Polishchuk and R. Vakil. Special thanks to 
 C. de Concini and A. Vistoli.

\section{Notation and terminology}
In this paper, the term {\em complex surface} refers to a smooth, connected, complex-analytic
two dimensional manifold with countable topology.

- $D'^p$ the sheaf of $p$-currents on a smooth manifold $M$; if the topology
is second countable, then the ``de Rham" complex
$\comp_M \to D'^{\bullet}$ is a fine resolution of $\comp_M$. This is the only reason why we require the topology  of $X$ to be countable.

- $\D := \{\, (z_1,z_2)\subset \comp^2 \, , \; |z_i|<1, \; i=1,\,2  \,\}$, the unit  a bi-disk.

- $X$ a complex  surface. 

- $\SS_n$ the symmetric group over $n$ elements.

- $P(n)$ the partitions of $n$.
 We use two standard pieces of notation.

\noindent 
$\bnu:=(\nu_1, \ldots, \nu_k)$, with $\nu_j>0$ and $\sum_{j=1}^k{\nu_j}=n$.

\n
The same partition can be represented as follows

\noindent
(a-notation) $\ba{a}=\ba{a}(\bnu) =(a_1, \ldots, a_n)$, where $a_i$ is the number of times that the number $i$ appears in the partition
$\bnu$. Note that $\sum_{i=1}^n{ia_i=n}$.

- $\lambda (\bnu) = \lambda (\ba{a})= k = \sum_{i=1}^n{a_i}$ the length of the partition.

- $X^n$, $X^{(n)}$ and $X^{[n]}$ the cartesian product, the symmetric product and the Douady space of zero-dimensional analytic subspaces of $X$
 of length $n$.

- $\pi = \pi_n: X^{[n]} \to X^{(n)}$, the natural Douady-Barlet morphism.

- $X^{(n)}_{(\bnu)}$ or $X^{(n)}_{(\ba{a})}$ the locally closed smooth
subspaces
of $X^{(n)}$ locus of points of the form $\sum_{j=1}^k \nu_ix_j$,
where the $x_j$ are pairwise distinct.

- $X^{[n]}_{(\bnu)}$ (or $X^{[n]}_{ (\ba{a}) }$) $:=$ 
$\pi^{-1} ( X^{(n)}_{(\bnu)} )$, where the pre-image is taken with
the induced reduced structure. 

- $\D^{[n]}_o: = \pi^{-1}(no)$, the closed reduced analytic subspace of $\D^{[n]}$ locus of subspaces of $\D$ of length $n$ supported at the origin $o\in \D$ (the so-called {\em punctual Hilbert schemes}).

- $X^{(\ba{a})}$, $\ba{a} \in P(n)$, the spaces
$\prod_{i=1}^n X^{(a_i)}$.

- $K_{\ba{a}}: X^{(\ba{a})} \to   \overline{X^{(n)}_{(\ba{a})}}$
the natural map sending
$$
( x_1^1 + \ldots + x_{a_1}^1, \, \ldots \, , x_1^n + \ldots + x_{a_n}^n)
\longrightarrow
\sum_{i=1}^n{i( x_1^i + \ldots + x_{a_i}^i)}.
$$ 

-
$$
X^{(n)}_l := \overline{
\coprod_{\lambda(\ba{a})= l} X^{(n)}_{ (\ba{a}) }  }
 =
\coprod_{ \lambda(\ba{a})= l} \overline{ X^{(n)}_{ (\ba{a}) }   } 
=\coprod_{\lambda(\ba{b})\leq l} X^{(n)}_{ (\ba{b}) }
$$

- $K_l= \coprod_{\lambda(\ba{a})=l} K_{\ba{a}} \longrightarrow X^{(n)}_l$.

\section{The Douady space of a bi-disk}
\label{dsbd}
 
In this section we recall the definition and 
main properties of the Douady space
$X^{[n]}$ of ``0-dimensional subspaces of length n" of a complex surface $X$;
see \S\ref{douadyxn}. 
We give  a self-contained construction, i.e. without assuming the main existence
result of Doaudy \ci{do},  of the Douady space and of the associated Douady-Barlet morphism in the case of  ${\comp}^2$
and of the bi-disk $\Delta$;
see \S\ref{toymodel}. We  do the same thing for {\em any}
 complex surface using a patching argument; see \S\ref{patch}.

\subsection{The Douady spaces $D(X)$ and $X^{[n]}$}
\label{douadyxn}

\bigskip
Let $An$ be the category of analytic spaces and
$X$ be an analytic space. A {\em family of compact subspaces of $X$ parameterized
by an analytic space $S$} is, by definition, 
an analytic subspace of $S\times X$ proper and flat over $S$.
Let $\Phi$ be the contravariant functor  which assigns to any analytic
space $T$ the set of families of compact subspaces of $X$ parameterized by $T$.
A fundamental result of A. Douady \ci{do} asserts that the functor
$\Phi$ is representable, i.e. there exists, unique up to unique isomorphism, a complex space $D(X)$ such that  $\Phi \simeq Mor_{An}(-, D(X))$ as functors;
the analytic space $D(X)$  is called
the {\em Douady space of $X$}. By definition, there is a universal flat
family
$u: {\cal U}_X \longrightarrow D(X)$.

\medskip
Let $X$ be a quasi-projective complex scheme. Let $H(X)$ be the associated Hilbert scheme and $H(X)_{An}$ be the associated analytic space.
The universal
family over $H(X)$ gives a holomorphic family over 
$H(X)_{An}$ and therefore, by virtue of the universal property and of
Chow Theorem,  there is a natural  bijective
morphism   $f: H(X)_{An} \to D(X)$.
First order algebraic deformations induce analytic ones so that
 the differential
$df$ is  injective at every point. Let $t$ be a point in $H(X)_{An}$ and
$f(t)$ be the corresponding point in $D(X)$. By ``GAGA"
principles the corresponding Zariski tangent spaces have the same dimension.  It follows  that $df$ is also surjective at every point, i.e. $df$ is an isomorphism at every point.
In particular,  $f: H({\Bbb A}^2_{\comp})_{An}
\simeq D(\comp^2)$: in fact these spaces are both smooth (see Theorem \ref{forse}). The same is true for every smooth algebraic surface, as long as we consider
those components parametrizing zero-dimensional components.
We thank R. Vakil for a useful conversation concerning this point.

\begin{rmk}
\label{pathology} 
{\rm 
If $X$ is projective, then every connected component
of $H(X)$ is projective and, in particular, it is  compact.
This is no longer true for an arbitrary (e.g. non K\"ahler)
 compact complex manifold; see \ci{ue}.
%the example is a compact threefold $X$ not in the Fujiki class $\cal C$; 
%the Douady space $D(X)$ admits a  smooth irreducible non compact  component; %this component
%parameterizes integral curves. 
But things can get even worse.
Consider Douady's example of the complex singular space $X$ constructed by identifying
a line and a conic in $\pn{3}$. One can construct
a connected component $D$  of $D(X)$ with an infinite
number of  irreducible
components $D_i$; each $D_i$ is  compact, but $\sup_i{ \{\dim{D_i}\} }=+\infty$
and the number  
of   connected components of the members of the family over $D$  can be made to be arbitrarily large.
This example grew out of conversations with  T. Graber, T. de Jong and
 R. Vakil.
The usual example, due to Hironaka, of a compact smooth threefold which is not a scheme
leads to pathologies for the Barlet space of cycles, but, apparently, not
for the Douady space.
}
\end{rmk}

The kind of pathologies described above  do not occur
for the Douady space of  zero-dimensional subspaces of a complex surface
(and more generally of a smooth $n$-fold).

\n 
Let $X$ be an analytic space and consider  zero-dimensional analytic subspaces: the flatness of $u: {\cal U}_X \to D(X)$ decomposes the open and closed subset of $D(X)$
corresponding to zero-dimensional families into the disjoint union
of the connected components over which the family has a fixed degree.
Let $\Phi^n$ be the sub-functor of  $\Phi$ corresponding to those
families which are flat and finite of degree $n$ over the base. By what above,
the functor $\Phi^n$
is represented by an open and closed subspace $X^{[n]} \subseteq D(X)$.
If $X$ is a compact complex surface, 
then $X^{[n]}$ is compact; see Theorem \ref{dcm}.

%then it is easy to show that the
%Douady-Barlet morphism  $\pi: X^{[n]}_{red}=X^{[n]} \to X^{(n)}$ (cfr. %\ref{dcm} and \ref{iversen})
%is projective.  In particular, if in addition $X$ is compact, then
%so is each    space  $X^{[n]}$.
%If $X$ is smooth of higher dimension, then
%the morphism  $\pi : X^{[n]}_{red} \to X^{(n)}$ is  proper and even projective, %locally over the
%base.

\medskip
Let $j: U \to X$ be an open immersion  of complex spaces. Then $D(U)$ sits
naturally (with respect to $j$) inside of $D(X)$. In particular,
$\D^{[n]} \subseteq  {\comp^2}^{[n]}= ({{\Bbb A}_{\comp}^2})^{[n]}_{An}$. This open immersion is  made explicit
below via the ``toy model" which we  describe momentarily.

\subsection{The toy model for ${\comp^2}^{[n]}$ and $\D^{[n]}$}
\label{toymodel}
An explicit construction of the Hilbert schemes of $n$-points $({\Bbb A}^2_{\comp})^{[n]}$, based on its existence, which was proved by A. Grothendieck, can be found in \ci{naknotes}
\S1. 

We now show how the construction satisfies the universal property so that
we provide a self-contained and complete construction of $({\Bbb A}^2_{\comp})^{[n]}$ which is independent of the usual general 
existence result. The proof works algebraically as well as analytically
and gives  the existence of the Douady spaces 
${\comp^2}^{[n]}$ which is independent of \ci{do}. 
This construction also establishes the existence of the Douady
spaces $\D^{[n]}$  and it identifies them concretely as sitting inside of ${\comp^2}^{[n]}$. 

\medskip
The construction of \ci{naknotes} \S1 proves  that 
$$
{\cal P}^n: = \{ (A,B,v) \in 
{\bf gl}(n) \times {\bf gl}(n) \times {\comp}^n \; \; | \; \; [A,B]=0,\, 
\comp^n= Span \{  A^kB^lv  \}_{k,l\, \in \nat}  \}
$$ 
is connected and smooth,  that it carries a flat family 
$w: {\cal W}^n \to {\cal P}^n$
of degree
$n$, that ${\bf GL}(n)$
acts freely on ${\cal P}^n$ by $G(A,B,v)=(GAG^{-1},GBG^{-1},Gv)$ and 
that  the quotient ${\frak g}: {\cal P}^n \to {\cal Q}^n:=
{\cal P}^n/{\bf GL}(n)$ :  1) exists, 2) is connected and 
  smooth and  3) that the family on ${\cal P}^n$ goes to the quotient
and defines
 a flat  family $v: {\cal V}^n \to {\cal Q}^n$ of degree $n$. We now show that the family over ${\cal Q}^n$ is universal.  We prove the result in $An$.
\begin{tm}
\label{forse}
Let $X=\comp^2$. The functor $\Phi^n$ is represented by $({\cal Q}^n_{An},
v)$, i.e. ${\comp^2}^{[n]}$ exists and is isomorphic to ${\cal Q}^n$.
\end{tm}
{\em Proof.}
Since the functor ${\Phi}^n$ is a sheaf of sets with respect to the classical topology, it is enough to show that ${\cal Q}^n$ enjoys the universal property
with respect to {\em germs of families} $\varphi: F \to S$ parametrized
by analytic germs $(S,s)$, i.e.: we need to show that if $F\subseteq S \times \comp^2$ is finite and flat of degree $n$ with respect to the first projection, then there exists a unique morphism of germs $\alpha(\varphi): (S,s) \to ({\cal Q}^n, q)$ such that
the germ of families $\alpha(\varphi)^* {\cal V}^n$ {\em is} $\varphi:
F \to S$.

\n
Let $\varphi: F\to S$ be such a germ.
The coherent $\odix{S}$-module
$M:= \varphi_*\odix{F}$ is free of rank $n$.

\n
The coordinates $(z_1,z_2)$, acting by multiplication on $M$,    define
two commuting $\odix{S}$-linear endomorphism $T_1$ and $T_2$ of $M$.
The distinguished section $1 \in M$ generates $M$ under the action 
of the momomials $T_1^kT^l_2$.

\n
Choose a $\odix{S}$-linear isomorphism $\e: M \to \odix{S}^{\oplus n}$.
By specializing at $s$, and using the trivialization $\e$, we get
$p:=(A,B,v)\in {\cal P}^n$. By virtue of \ci{fischer} 0.21, this datum is equivalent to giving a morphism of germs $\tilde{\alpha}(\varphi)_{\e}: (S,s) \to (
{\cal P}^n, p)$. By construction $\tilde{\alpha}(\varphi)_{\e}^*
{\cal W}^n$ {\em is} $\varphi: F \to S$.

\n
By composition we get ${\alpha}(\varphi):= {\frak g}\circ \tilde{\alpha}(\varphi)_{\e}: (S,s) \to ({\cal Q}^n, q:={\frak g}(p))$
with the property that
${\alpha}(\varphi)^*
{\cal V}^n$ {\em is} $\varphi: F \to S$. The morphism 
${\alpha}(\varphi)_{\e}$ is independent of $\e$ so that we denote it by
${\alpha}(\varphi)$. 
\n
To prove uniqueness, note that if $\eta : (S,s') \to ({\cal Q}^n,q')$ is any morphism of germs, then we get a germ of  families
 $\s: \eta^*{\cal V}^n \to S$,
and $\alpha(\s)=\eta$. \blacksquare

\bigskip
We now recall the definition of the Douady-Barlet map.
The starting point is a ``doubled'' version of the map
$ {\bf gl}(n) \to  {\bf gl}(n)/Ad \simeq {\comp}^{(n)}$
which associates with a matrix its  characteristic polynomial.
By a classical theorem in Invariant Theory
(cfr. \ci{Weyl}), the ring of invariant functions of $2n$ variables,
 $((x_1,y_1),\ldots ,(x_n,y_n))$, by the action of the symmetric group
$\SS_n$ is generated by the functions
$f_{k,l}( (x_1,y_1),\ldots , (x_n,y_n))=\sum_i x_i^ky_i^l$.\,
In other words, 
${{\comp}^2}^{(n)}=
Spec_{An} \, {\comp}\,[\, x_1,y_1,\ldots , x_n,y_n \,]^{\SS_n}=
Spec_{An} \, {\comp} \, [\,  \ldots , f_{k,l}  ,\ldots \,]$.

Let  $(A,B,v)\in {\cal P}^n$. Since $A$ and $B$ commute,
we can let the group ${\bf GL}(n)$ act so that they are both 
in triangular form. If $x_1,\ldots ,x_n$ and $y_1,\ldots  , y_n$ are 
the respective diagonal terms, well defined up to a permutation,
then $f_{k,l}( x_1,y_1,\ldots , x_n,y_n)=\sum_i x_i^ky_i^l=
Trace\, (A^k B^l)$. These functions  are   ${\bf GL}(n)$-invariant. 
We can therefore define  a   
morphism $\pi: {{\comp}^2}^{[n]} \lorw {{\comp}^2}^{(n)}$,
which is called the {\em Douady-Barlet} 
morphism, or the Hilbert-Chow morphism in the algebraic category.
In set-theoretic terms this map associates with a subspace of length $n$
of $ { \comp }^2$ its support, counting multiplicities, seen as a point in   ${{\comp}^2}^{(n)}$.
The Douady-Barlet morphism is easily seen to be proper and its fibers are connected by Zariski Main Theorem.

\medskip
Let  $U$ be an open subset of ${\comp}^2$. Since the quotient map
$({\comp}^2)^n \lorw {{\comp}^2}^{(n)}$ is finite and therefore open,
$U^{(n)}$ is open in ${{\comp}^2}^{(n)}$  and 
$U^{[n]}$ ($\pi: U^{[n]} \to U^{(n)}$, respectively)
can be naturally identified with $\pi^{-1}(U^{(n)})$ 
($  \pi: \pi^{-1} (U^{(n)}) \to U^{(n)}$, respectively).
In particular, we have

\begin{pr}
\label{dn}
 Let $\Delta=\{(z_1,z_2)\in {\comp}^2: |z_i| < 1 , \, i=1,\,2 \}$ 
 be the two-dimensional bi-disk. The Douady space  $\Delta^{[n]}$
 can be described as the quotient by the action of ${\bf GL}(n)$ of 
$$
{\cal P}^n_{\Delta}: = \left\{ (A,B,v) \in {\cal U}^n  \; | \;  
 \hbox{ the 
eigenvalues of }A \hbox{ and } B \hbox{ have modulus smaller than one} \right\}.
$$ 
\end{pr}

\begin{rmk}{\rm The natural action of ${{\comp}^*}^2$ on ${\comp}^2$
induces an action on ${{\comp}^2}^{[n]}$ given by
$(\lambda_1,\lambda_2)(A,B,v)= (\lambda_1 A,\lambda_2 B,v)$.
The subspace $\Delta^{[n]}$ is $(S^1)^2$-invariant under this action.
The fixed points are easily determined and they correspond to
the partitions of $n$;
see \ci{e-s} \S3 and \ci{naknotes} \S5.}
\end{rmk}

\bigskip
The following is a nice consequence of the toy model developed above.

\begin{pr}
\label{ret}
The manifolds ${\comp^2}^{[n]}$ and $\D^{[n]}$ are homeomorphic to each other.

\noindent
In particular, ${\comp^2}^{[n]}$ and $\D^{[n]}$ have the 
same Betti numbers. 
\end{pr}
{\em Proof.} Let $\phi(A,B ) := \hbox{maximum of the absolute values of the eigenvalues of }
A \hbox{ and }B$. 
Consider the application
$\Phi :  {\Delta}^{[n]} \longrightarrow {{\comp}^2}^{[n]}$ defined
as follows:
$$
\Phi (A,B,v)=\left(\frac{1}{1-\phi(A,B)}A\,, \, \frac{1}{1-\phi(A,B)}B\,
, \, v \right).
$$
This map is ${\bf GL}(n)$-invariant and defines a homeomorphism.
\blacksquare

\medskip
At this point we know the Betti numbers of $\D^{[n]}$ by virtue
of Proposition \ref{ret} 
and of \ci{e-s} Theorem 1.1 (since loc. cit. has a typo, the reader is referred to
\ci{naknotes} 5.9). However, our analysis of $\pi: X^{[n]} \to X^{(n)}$ will require a  certain basis  for the cohomology which reflects
the geometry of the stratified Douady-Barlet morphism. We shall obtain such a basis
via Nakajima's construction.  Compare with Remark \ref{comparebases}.

\subsection{Construction of $X^{[n]}$ and $\pi : X^{[n]} \to X^{(n)}$}
\label{patch}
We can now sketch a proof 
of the existence of $X^{[n]}$ and of $\pi: X^{[n]} \to X^{(n)}$
by a patching argument. In the algebraic case (and assuming the existence
of the Hilbert scheme) the corresponding assertions are due to
Fogarty \ci{fo}.
\begin{tm}
\label{dcm}
Let $X$ be a  complex surface. For every $n \in \nat$, the Douady space  
$X^{[n]}$ and the Douady-Barlet morphism $\pi : X^{[n]} \to X^{(n)}$
exist and  can be constructed by patching
using $\D^{[m]}$ and  $\pi: \D^{[m]} \to \D^{(m)}$, $m\leq n$. The space
$X^{[n]}$ 
is a connected $2n$-dimensional complex manifold. 

\n
The Douady-Barlet morphism
is projective.   In particular, if
$X$ is compact, then $X^{[n]}$ is compact.
\end{tm}
{\em Sketch of proof.} We freely use
the local descriptions given in Sction \ref{douadybarlet}. 
Let $z \in X^{(n)}$ be any point. There is a unique partition
$\bnu \in P(n)$ such that $z=\sum_{j=1}^{k=\lambda (\bnu)}{\nu_j x_j}
\in X^{(n)}_{(\bnu)}$, where the points $x_j\in X$ are pairwise distinct. Consider the complex manifold
$\prod_{j=1}^k{_{x_j}\D^{[\nu_j]}}$. This way, we obtain a set of charts which glue coherently  by the
universal property of  the Douady space (which we have already constructed for
$\D$) and define a connected  complex manifold $W^n$. 
 Similarly, one can check that the local Douady-Barlet morphisms glue
 and define a global proper map with  connected fibers $\pi : W^n \lorw X^{(n)}$.
By construction $W^n$
carries a family. It is elementary to check the universal property of this family using the corresponding fact for $\D$ and by the fact that 
the functor $\Phi$ is a sheaf for the classical topology. 

\n
The Douady-Barlet morphism is bimeromorphic. It is also projective, locally over $X^{(n)}$.
We conclude by the fact that there is a unique irreducible exceptional divisor, i.e.
the pre-image of the big diagonal in $X^{(n)}$. 
\blacksquare 

\begin{rmk}
\label{iversen}
{\rm An alternative definition of $\pi$ was given   
by Iversen \ci{iv}. More generally, D. Barlet
has constructed the so-called Barlet spaces of cycles $B(X)$ associated with any analytic space $X$ and a natural morphism $D(X)_{red} \to B(X)$ which, in our case, gives  $\pi: X^{[n]} \to X^{(n)}$. This explains the name 
Douady-Barlet for the morphisms $\pi$. 
One can show that this morphism is proper
if we consider zero-dimensional subspaces of length $n$ on a smooth manifold.
}
\end{rmk}

\subsection{The operators of Nakajima} 
\label{opnak}
Let $X$ be a complex surface. In what follows we work with rational 
singular cohomology. We now recall, for the convenience of the reader, Nakajima's construction
\ci{nak}, \ci{naknotes}
of the correspondences which realize geometrically an action
of the Heisenberg/Clifford super-algebra on 
$$
{\Bbb H}(X):= \bigoplus_{r\geq 0}{
H^*(X^{[r]})}.
$$
Let $Y$ be a topological manifold of dimension $m$; we denote rational 
cohomology simply by $H^*(Y)$ and we identify it freely with Borel-Moore homology $H^{lf}_{m-*}(Y)$ via Poincar\'e Duality.  

\medskip
For every  $r\geq 0$ and every $k>0$, define
a  closed and reduced  analytic subspace 
$T_k \subseteq X \times \H{r} \times \H{r+k}$
by setting
$$
T_k= \left\{ (x,\zeta_1,\zeta_2) : \, {\cal I}_{\zeta_2} \subseteq 
{\cal I}_{\zeta_1},\, {\cal I}_{\zeta_1}/{\cal I}_{\zeta_2} \hbox { is supported at } x \right\},
$$
where ${\cal I}_{\zeta_i}$, $i=1,2$, are  the coherent sheaves of ideals associated with
the closed analytic space $\zeta_i \subseteq X$ $i=1,2$.

\n
By \ci{naknotes} \S8.3, $T_k$ has a unique irreducible reduced component
$Z_k$ of (expected) dimension $2 + 2r + (k-1)$ and all other irreducible components, if there are any at all, have strictly lower dimension.
We ignore whether $T_k$ is irreducible (i.e. $T_k=Z_k$) or not. However, 
this does not cause any problem in what follows.

\medskip
We are implicitely making use of Brian\c{c}on's irreducibility result
\ci{br}, V.3.3 as in \ci{nak} or \ci{naknotes}. However, this turns out the be not necessary; see Section \ref{punctual}. At this stage, we prefer to
use this result for clarity of exposition.

\bigskip

The projections  $p_{12}: X \times \H{r} \times \H{r+k}
\lorw X \times \H{r} $
and 
$p_3:X \times \H{r} \times \H{r+k}
\lorw  \H{r+k}$
are proper when restricted to $Z_k$. Therefore we can define the two correspondences, 
one transpose of the other,
$$ 
{\cal P}[k]:= p_{3*}(p_{12}^*(...) \cap Z_k):
             H^*(X \times  \H{r}) \lorw H^*(\H{r+k}),
$$
$$ 
{\cal R}[k]:= p_{12*}(p_{3}^*(...) \cap Z_k):
             H^*(\H{r+k}) \lorw H^*(X \times \H{r}),
$$
or, equivalently, 
$$ 
\tilde{\cal P}[k]:
             H^*(X) \lorw Hom(H^*(\H{r}),H^*(\H{r+k})),
$$
$$
\tilde {\cal R}[k]= 
             H^*_c(X) \lorw Hom(H^*(\H{r+k},H^*(\H{r})).
$$
Putting all $r$'s together we obtain  endomorphisms of ${\Bbb H}(X)$.
The endomorphisms corresponding to 
$\alpha \in H^*(X)$ and $\beta \in H^*_c(X)$ will be denoted by
$\tilde{\cal P}_{\alpha}[k]$ 
and 
$\tilde {\cal R}_{\beta}[k]$, respectively.

Since $\dim{Z_k}=k+1+2r$ and $\dim{  X \times \H{r} \times \H{r+k}} =2(k+1+2r)$,
it is easy to compute that ${\cal P}[k]$ increase the degree by $2(k-1)$
and ${\cal R}[k]$ decrease the degree by the same amount, and 
similarly if we 
look at the correspondences on homology groups or on locally finite 
homology groups. 

\medskip
The following is an explicit, simple, but important examplification of how
these operators work.

\begin{lm}
\label{afact}
Let $X$ be any complex surface,
 $\bnu= (\nu_1, \ldots, \nu_k) = \ba{a} \in P(n)$, $\ba{a}!:=
\prod_{i=1}^n{a_i!}$,  $[\overline{X^{[n]}_{(\bnu)}}]
\in H^{2(n-k )}( X^{[n]},\rat)$
be the fundamental class of $\overline{X^{[n]}_{(\bnu)}}$. Then
$$
\tilde{\cal P}_{[X]}[\nu_1] \circ \ldots  \circ \tilde{\cal P}_{[X]}[\nu_k] ({\bf 1})=
\ba{a}! \,[\overline{X^{[n]}_{(\bnu)}}].
$$
\end{lm}
{\em Proof.} By induction on the length $k$ of $\bnu$, it is enough to prove that 
$$
\tilde{\cal P}_{[X]}[\nu_k]\, ( [ \overline{X^{ [n- \nu_k] }_{ (\nu_1, \ldots, \nu_{k-1}) }} ] ) = a_{\nu_k} \,[\overline{X^{[n]}_{(\bnu)}}].
$$
To compute the left hand-side we first consider the intersection
$T:= p_1^{-1}(X) \cap p_2^{-1}( \overline{X^{ [n- \nu_k] }_{ (\nu_1, \ldots, \nu_{k-1}) }})
\cap Z_k$ inside of $X \times X^{[n-\nu_k]} \times X^{[n]}$. The space $T$ is irreducible of the expected dimension and
the intersection is transverse at a general point of $T$.
The third projection ${p_3}_{|}: T \to X^{[n]}$ is a proper, generically finite morphism onto its
image $\overline{X^{[n]}_{(\bnu)}}$.  We can determine the generic degree
at a general point. This  degree  is $a_{\nu_k}$.
\blacksquare

 \begin{rmk}
\label{hodge1}
{\rm Let  $X$ be a Zariski-dense open subset of a compact K\"ahler surface $Y$, e.g. X is quasi-projective. Then $X^{[n]}$ is a Zariski-dense open subset of
the compact K\"ahler manifold $Y^{[n]}$.
It follows that
$H^*(X)$ and $H^*(\H{r})$ have a Hodge decomposition.
Since the class $Z_k$ has type $(k+2r+1,k+2r+1)$ and 
$p_{3*}$ has degree $(-(2+2r),-(2+2r))$, it follows that
${\cal P}[k]$  
and 
${\cal R}[k]$, 
have degree $(k-1,k-1)$ and 
$(1-k,1-k)$ respectively.

\n
If $\alpha$ has type $(p,q)$, then 
$\tilde{\cal P}_{\alpha}[k]:
H^{*,*}(\H{r}) \lorw H^{*+p+k-1,*+q+k-1}(\H{r+k})$.} 
\end {rmk}

\bigskip
As anticipated,  Nakajima's construction and G\"ottsche's 
formula can be conveniently encoded in representation-theoretic terms.
For this it is necessary to introduce the language of 
super (=${\zed}_2$-graded)  algebras. A general reference is \ci{kac}.

\n
Let us consider $H^*(X)=H^{even} \oplus H^{odd}$ as a super vector space,
and define the Heisenberg superalgebra 
$$
{\cal H}(X) : = {\rat <c>} \oplus {\bf s}_{> 0}\oplus {\bf s}_{< 0}\, ,
$$
where
$$
{\bf s}_{> 0} =  \bigoplus_{i>0}H^*(X)[i], 
\qquad  \qquad
{\bf s}_{< 0} =\bigoplus_{i<0}H^*_c(X)[i],
$$
(here $[i]$ is just a place holder and does not refer to 
any shift) the $\zed_2$-grading comes from the one on $H^*(X)$, the $\zed$-degree of $c$ is zero, the $\zed_2$-degree of
$c$ is even, and  the defining properties are

\smallskip
a) $c$ is central,

b) ${\bf s}_{> 0}$ and ${\bf s}_{< 0}$ are supercommutative, 
i.e. the supercommutator is identically zero,

c) $\left[\, \beta[j],   \alpha[i] \, \right] \, =\,
\delta_{i,-j}  \, (-1)^{i-1}i \,  <\alpha,\beta> \,c$,
\smallskip

\n
where $<-,->$ is the canonical Poincar\'e pairing between cohomology 
and   cohomology with compact supports. Note that the factors
$(-1)^{i-1}i$ are not necessary and could be either omitted, or replaced
by another compatible system of factors. We have used them in view of
Theorem  \ref{supercommuting} where they appear naturally
in connection with the computation of a certain intersection number
which has been  determined in \ci{e-s2}.

\n
Since ${\bf s}_{>0}$ is supercommutative, the enveloping superalgebra
$U({\bf s}_{> 0})$ (which is constructed
like its non-super analogue by factoring the tensor algebra by the ideal
$x\otimes y -(-1)^{degxdegy}y\otimes x -[x,y]$, where $[-,-]$ is 
the supercommutator)
is the free commutative superalgebra on 
$H^*(X)$ isomorphic to 
$S^*({\bf s}_{> 0}^{even})\otimes \bigwedge^*({\bf s}_{> 0}^{odd})$. 
By its own  definition, ${\bf s}_{> 0} $
is  
${\nat}^+$-graded, and an ${\nat}$-grading is inherited by $U({\bf s}_{> 0})$.

If $\dim{H^*(X)}< \infty$, then a computation shows that
\begin{equation}
\label{qdegree}
 \sum_{n=0}^{\infty} q^n \dim{  U({\bf s}_{> 0})_n  } 
 =
\prod_{m=1}^{\infty}{
\frac{ (1+  q^m)^{dimH^{odd}(X)} \; }{
(1-  q^m)^{dimH^{even}(X)}\; }
}.  
\end{equation}

Denote by $[- , -]$ the supercommutator in
$ \hbox{\rm End} \, {\Bbb H}(X)$. The main theorem in \ci{nak}, or
\ci{naknotes}, \S8 is 

\begin{tm}
\label{supercommuting}
 The following commutation relations hold:
$$
[\, \tilde{\cal P}_{\alpha}[k],\tilde{\cal P}_{\gamma}[l]\,]=0,
$$
$$
[\, \tilde{\cal R}_{\beta}[k],\tilde{\cal R}_{\delta}[l] \, ]=0,
$$
$$[\, \tilde{\cal R}_{\beta}[l] , \tilde{\cal P}_{\alpha}[k] \, ]
            =\delta_{k,l}(-1)^{k-1}k<\alpha,\beta> \, Id_{{\Bbb H}(X)}.
$$ 
\end{tm}
Note that the ``reversed" order  in the third relation 
is due to the fact that we are using operators which are  dual to
 the ones that Nakajima defines in homology (cfr. \ci{naknotes}, \S8). Equivalently, the assignments:

\smallskip
\n
$\alpha[i] \lorw \tilde{\cal P}_{\alpha}[i]$ if  $i>0$ and $\alpha \in H^*(X)$,

\smallskip
\n
$\beta[i] \lorw \tilde{\cal R}_{\beta}[-i]$ if  $i<0$ and $\beta \in H^*_c(X)$,

\smallskip
\n
$c \lorw Id_{{\Bbb H}(X)}$ 

\smallskip
\n
exhibit ${\Bbb H} (X)$ as a representation of ${\cal H}(X)$.
For every $i>0$, the vector ${\bf 1} \in H^0(\H{0})$ is annihilated by  $\tilde{\cal R}_{\beta}[i]$ and is therefore
a highest weight vector.
Furthermore, it is immediate that, 
$\forall v \in {\Bbb H}(X)$, 
         $\tilde{\cal R}_{\beta}[i]\, v=0$ for 
                     $i$ big enough depending on $v$.

%\medskip
%Under these assumptions it can be shown  (cfr. \ci{kac}
%\S 9.13) that ${\Bbb H}(X)$
%is a direct sum of representations isomorphic to 
%$U({\bf s}_{> 0})$ up to some shift.

\n
Denote by $U'$ the ${\cal H}(X)$-submodule  of ${\Bbb H}(X)$ 
generated by $\bf 1$.
It is well known that $U'$ is irreducible and isomorphic to
$  U({\bf s}_{> 0})$. 
Because of formula  (\ref{qdegree}),
G\"ottsche's formula
(see Theorem \ref{go}.(\ref{gofock})) becomes the statement 
$U'={\Bbb H}(X)$.

\medskip
Let us summarize Nakajima's construction and its link to G\"ottsche Formula.

\begin{tm}
\label{nakgo}
Let $X$ be a complex surface. Then 

\smallskip
\n
(\ref{nakgo}.1) ${\Bbb H}(X)$ is a highest weight representation
of ${\cal H}(X)$, geometrically realized by the operators $\tilde{\cal P}_{\alpha}[k]$
and 
$\tilde{\cal R}_{\beta}[l]$.

\smallskip
\n
(\ref{nakgo}.2) the representation is irreducible iff G\"ottsche formula holds for
$X$.
\end{tm}

\subsection{Irreducibility of the action  for $X=\D$ and a canonical basis for $H^*(\D^{[n]}, \rat)$}
\label{canbasis}

\bigskip

We now examine in detail the picture in homology in 
the case of 
$\comp^2$ and $\Delta$. They can be dealt with simultaneously.
Let $X=\comp^2$ or $X=\D$.

%\medskip
%A generator of $H_*(X)$ is the class 
%$[x]\in H_0(X)$ of  a point $x\in X$.
%The space  $\H{0}$ is a single point corresponding to the empty scheme; let %${\bf 1}$ denote the canonical generator of  
%$H_0((X)^{[0]})$. There is only one operator 
%${\cal P}[k]:=\tilde {\cal P}_{[x]}[k] \in  \mbox{End} \,{\Bbb H}(X)$.

%Similarly, the generator of  the dual space $H_*^{lf}(X)$ is the class
%$[X] \in H_4^{lf}(X)$ of $X$, and we have the corresponding operator 
%$ {\cal R}[k]:=\tilde {\cal R}_{[X]}[k] \in \mbox{End} \, {\Bbb H}(X)$.

%Let $\tilde \nu=\nu_1,\ldots, \nu_{\lambda(\tilde \nu)}$, with $ \sum \nu_j=n$,
%be a partition of $n$. By the very definition 
%of the correspondences $Z_k$,
 %the homology class 
%$P[\tilde \nu]\in H_{2(n-l(\tilde \nu))}(\H{n})$ defined as
%$ P[\tilde \nu]:=   P[\nu_1] \circ ...\circ P[\nu_{\lambda (\tilde \nu)}](
%{\bf 1})$
 %can be represented by the cycle
%$\pi^{-1}(\nu_1 x_1+ \ldots +\nu_{\lambda (\tilde \nu)}x_{\lambda(\tilde %\nu)})$, i.e. the fibre over the corresponding point, given the reduced %structure.

\bigskip
The Heisenberg algebra ${\cal H}(X)$ is in this case the standard one and the 
fundamental representation  $U({\bf s}_{> 0})$ is isomorphic 
to the space of polynomials in  infinitely many variables $p_i$
to which we assign degree $i$. See \ci{naknotes} \S8.
It follows therefore that, for every $n$,  $\dim{ U({\bf s}_{> 0})_n}=p(n)$, 
the number of partitions of $n$.

\medskip
We can now give a quick and self-contained  
proof of the formula for the Betti numbers of the Douady space of 
$\comp^2$ or $\Delta$, very much in the spirit of this paper.
This formula, and especially the construction of the explicit basis  
derived from it,
are heavily used in the sequel of this paper.
For the original proof of the formula for $X=\comp^2$,  see \ci{e-s} or  
the slightly  different  method
in \ci{naknotes}.

\begin{tm} 
\label{pobd}
Let    $X=\comp^2$ or $X=\Delta$. 
Then 

\smallskip
\n
(\ref{pobd}.1) ${\Bbb H}(X)$ is an irreducible
highest weight vector representation
of the Heisenberg algebra with  highest weight vector  the generator
${\bf 1} \in H^0(X^{[0]})\simeq \rat$.

\smallskip
\n 
(\ref{pobd}.2)
The Poincar\'e polynomials of $X^{[n]}$ and their generating function are, respectively: 
$$
P_t(X^{[n]}) = \sum_{\bnu \in P(n)}{t^{2n - 2\lambda (\bnu)}}, 
\qquad and \qquad 
 \sum_{n=0}^{\infty} P_t ( X^{[n]} ) \, q^n  
 =
\prod_{m=1}^{\infty}{
\frac{ 1 }{
1-  t^{2m-2}q^m}.
}  
$$

\smallskip
\n
(\ref{pobd}.3)  The cycles
$\overline{X_ {(\tilde \nu)}^{[n]}}\in H^{lf}_{2(n+\lambda(\tilde \nu))}(\H{n})
= H^{2(n-\lambda(\tilde \nu))}(\H{n})$ form a basis for $H^*(\H{n})$.

\end{tm}
{\em Proof.} Since the generator  of $H^0(X^{[0]})$ is a highest weight vector,
it follows that ${\Bbb H}(X)$ contains a subrepresentation isomorphic to
$U({\bf s}_{> 0})$. In order  to prove the irreducibility, it is 
enough to prove that
$\sum_i \dim{ H^i(\H{n}) } =\dim{ U({\bf s}_{> 0})_n } = p(n)$, for every $n$.
As we already noted, $({\comp}^2)^{[n]}$ has an action 
of $ ({\comp}^*)^2$ with isolated fixed points.
It follows that the odd Betti numbers vanish and that the 
sum of the even Betti numbers is equal to the number of these fixed points.
In this case they correspond to monomial ideals of colength $n$ 
and 
the number of these is exactly $p(n)$. This proves the first part also
for $\D$ by virtue of Proposition \ref{ret}.

\smallskip
\n
Taking into account the degree properties of the operators
$\tilde{\cal P}_{[X]}[l]$, the second statement 
follows immediately from the first one.

%just amounts to the observation made above that 
 %$P[\tilde \nu]\in H_{2(n-l(\tilde \nu))}(\H{n})$ and one counts easily these %classes.

\smallskip
\n
The third one follows from Lemma \ref{afact}.
\blacksquare
%To prove the third, observe that we have,    
%by the 
%result of Ellingsrud-Stromme \ci{e-s2}
%$ P[\tilde \nu]\cap M[\tilde \nu]=
%\prod_{i=1}^{l(\tilde \nu)}(-1)^{k_i-1}k_i$,
%whereas the other intersection numbers are zero because 
%the supports are disjoint. \blacksquare

\begin{rmk}
\label{comparebases}
{\rm 
Theorem \ref{pobd} gives a different proof of \ci{e-s} Theorem 1.1.(iii)
and determines a  basis of elements for ${\Bbb H}(\comp^2)$ different
from  the one
in \ci{e-s1}. 
The basis (\ref{pobd}.3) for ${\Bbb H}(\D)$ is used in an essential way in our study of the Douady-Barlet morphism.
}
\end{rmk}

\begin{rmk}
\label{nicefact}
{\rm
Let $X$ be either $\comp^2$, or $\D$. The algebra ${\cal H}(X)$ is the standard
Heisenberg algebra. G\"ottsche Formula holds for $X$. See \ci{e-s}
and Proposition \ref{ret}; see also Theorem \ref{pobd}.
It follows that ${\Bbb H}(X) = <{\bf 1}>_{{\cal H}(X)}$. Recall that 
$<{\bf 1}>_{{\cal H}(X)} \simeq \rat [t_1, t_2, \ldots ]$ as ${\cal H}(X)$-modules; see \ci{naknotes}, \$8 for example, and assign
$t_i$ to $\tilde{\cal P}_{[X]} [i] ({\bf 1})$. Let
$\bnu \in P(n)$. By virtue
of Lemma \ref{afact} and after an obvious normalization, we see that  the monomial $t_{\nu_1} \, \cdots \, t_{\nu_k}$ corresponds to the class of $\overline{X^{[n]}_{(\bnu)}}$.  This  fits nicely with Theorem 
\ref{pobd}.3.
}
\end{rmk}

\subsection{Punctual Hilbert schemes: Betti numbers and irreducibility}
\label{punctual}
Let $n \in \nat$, $p\in X$  and $X^{[n]}_{p}$ be the fiber of $\pi$ over the
point $n  p \in X^{(n)}_{(n)}$ with the induced reduced structure. These spaces are analytically 
isomorphic to
the corresponding spaces $\D_o^{[n]}$ when $X$ is $\D$ (or $\comp^2$) and $p$ is the origin
$o$.
In particular they are  projective. They are called {\em punctual Hilbert schemes}.

\medskip
We can give a new  proof, based on Theorem \ref{pobd} of the following basic results 
\ci{e-s}, Theorem 3.1.(iv) and \ci{br} V. 3.3. concerning some of the topology
of these spaces.

%Let us introduce some notation in 
%\ci{e-s}. Let $m$ and $q$ be non-negative integers, $P(m,q)$ be the number of %sequences
%$q\geq b_0 \geq b_1 \geq \ldots \geq b_m =0$ such that $\sum{b_j}=m$.
%Note that if $q\geq m$, then $P(m,q)=|P(m)|$. Let
%$P(m,q)=0$ if either $m$ or $q$ is negative.

\begin{cor}
\label{bettibriancon} (Cf. \ci{e-s})
The Poincar\'e polynomials of $\D^{[n]}_o$
and their generating function are as follows:
$$
P_t(\D^{[n]}_o) = \sum_{\bnu \in P(n)}{t^{2n - 2\lambda (\bnu)}}, 
\qquad and \qquad 
 \sum_{n=0}^{\infty} P_t ( \D^{[n]}_o ) \, q^n  
 =
\prod_{m=1}^{\infty}{
\frac{ 1 }{
1-  t^{2m-2}q^m}.
}  
$$
\end{cor}
{\em Proof.}
It follows at once  from cohomology and base change for proper maps
and Theorem \ref{pobd}.3. 
\blacksquare

\bigskip
The space $\D^{[n]}_o$ is irreducible  by a result of Brian\c{c}on; see \ci{br} V.3.3;
this result is re-proved by different methods in \ci{e-s}, Corollary 1.2.

We want to observe that Briancon's result can be re-proved, in the spirit of this paper,
using Theorem \ref{pobd} above, exactly as \ci{e-s}, Theorem 1.1.iv is used to prove
\ci{e-s}, Corollary 1.2. Theorem \ref{pobd}  is based on the use of Nakajima construction which uses Brian\c{c}on's result. (It should be pointed out
that the hard part of this irreducibility result has to do
with bounding from below certain loci; once that has been achieved,
the irreducibility follows by standard arguments  not involving
the indirect argument  reproduced below; besides Brian\c{c}on's original proof, we would like to mention work of Gaffney, Gaffney-Lazarsfeld, Granger and Iarrobino;
see \ci{ia3}.)
 To avoid biting our tail we observe that
a modification  in the definition of the varieties $Z_k$ in \S \ref{opnak}
 allows us to modify slightly  Nakajima's construction without using the irreducibility of
punctual Hilbert schemes.

Let $CX^{[n]}_p$ be the closure in $X^{[n]}_p$ of the locus of curvilinear subspaces
supported at $p$. It is irreducible of dimension $n-1$ (cfr. the end of \S1 in \ci{naknotes}).
Note that the irreducibility of $\D^{[n]}_o$ is equivalent to $CX^{[n]}_p=X^{[n]}_p$.

Define:
$$
S_k:= \{   (x,\zeta_1, \zeta_2) \in X \times X^{[n]} \times X^{[n+k]} \, | \; 
x \notin Supp(\zeta_1), \, \zeta_2 = \zeta_1 \coprod \zeta_2', \, \zeta_2' \in CX^{[k]}_p        \},
$$
and $Z'_k=\overline{S_k}$. Nakajima's construction works if we replace $Z_k$ by $Z_k'$.
We thus have

\begin{cor}
\label{brianc}
(Cf. \ci{br}, \ci{e-s}) The spaces $\D^{[n]}_o$ are irreducible.
\end{cor}
{\em Proof.}
By a result of Gaffney-Lazarsfeld \ci{ia2}, Theorem 2,  all the irreducible components
of $\D^{[n]}_o$ have dimension at least $n-1$. By virtue of Theorem
\ref{pobd},  we have that   (i)  $b_{2(n-1)}(\D^{[n]}_o) =1$
and (ii) $b_t( \D^{[n]}_o )=0$, for every $t> 2(n-1)$.
The irreducibility follows easily.
\blacksquare

\section{Study of the local structure of the 
Douady-Barlet morphism}
\label{douadybarlet}
This section consists of the detailed analysis of the stratified
morphism $\pi : X^{[n]} \to X^{(n)}$.

\n
The goal is to define the morphism of complexes $\Psi$;
see Proposition \ref{prpsi}.

\n
The morphism $\Psi$ is defined via its components $\Psi^{2h}$ (the odd ones are zero). 

\n
The morphisms of sheaves $\Psi^{2h}$ are defined  
by establishing natural identifications between combinatorial
and geometric objects arising from our analysis of the stratified morphism. In short,
we first fix a point $z \in X^{(n)}$ and small natural euclidean neighborhoods ${\frak U}$
of $z$, then we identify  the fibers over $z$  of {\em all} the normalization morphisms
$K_{\ba{a}}$, for those $\ba{a}$ such that $\lambda ( \ba{a})=n-h$, with the closed currents of integration
associated with  \ref{pobd}.3, which give a basis  for the cohomology
in degree $2h$
of the open sets $\pi^{-1}({\frak U})$.

\n
The Decomposition Theorem \ref{decel} will be the statement
that $\Psi$ is a quasi-isomorphism.

\subsection{Partial ordering on $P(n)$}
\label{partial}
Let $\bnu:=(\nu_1, \ldots, \nu_k)$, $\bnu':=
(\nu_1', \ldots, \nu_l')$ be two partitions of $n$.
We say that $\bnu' \geq \bnu$ if there exists a partition of the set
$\{\nu_1', \ldots, \nu_k' \}$ into $l$ disjoint  subsets
${\Gamma}_j:= \{\mu_1^j, \ldots, \mu_{t_j}^j \}$ such that
 the elements of  ${\Gamma}_j$ form  a partition of $\nu_j$, for
every $j=1, \ldots, k$.
Geometrically: $\bnu' \geq \bnu$ iff $X^{(n)}_{(\bnu)} \subseteq
\overline{X^{(n)}_{(\bnu')}}$.

\begin{ex}
{\rm 
 Let $n=6$; one sees immediately that  the pairs
$(3,1,1,1)$
and $(2,2,2)$, $(4,1,1)$ and $(3,3)$, $(2,2,2)$ and $(5,1)$, $(2,2,2)$
and $(3,3)$  are not comparable in the sense that neither element of 
the pair is greater than or equal 
to the remaining  element of the pair. 
}
\end{ex}

\noindent
It may be useful to arrange all partitions of $n$, in $\bnu$-notation,
 in columns left to right 
following the  decreasing length   and draw either solid
arrows or dotted segments from the elements of one column
to the ones of the column immediately to the right
according to whether or not two partitions are comparable.
One finds that $\bnu \geq \tilde{\mu}$ iff
there is a path given by solid arrows from left to right
connecting $\bnu$ and $\tilde{\mu}$.

\subsection{A fundamental system of neighborhoods on $X^{(n)}$}
\label{basics}
\noindent
Recall that, given a point $y$ in a topological space $Y$
and  an indexing set $I$,  a {\em
fundamental system of neighborhoods of $y$ in $Y$ labelled by $I$} 
is a collection $U_i$, $i\in I$
of open neighborhoods of $y$, 
such that
given any open set in $Y$ containing $y$, there exists an index $j \in I$
such that $U_j \subseteq U$. Given such a system for every point
$y$ in $Y$, the collection of these open sets forms a basis 
for the topology of $Y$.

\smallskip
\noindent
Let $\bnu \in P(n)$ and $z=\sum_{j=1}^k{\nu_j x_j} \in X^{(n)}_{(\bnu)}$, i.e.
the points $x_j \in X$ are pairwise distinct. 
Note that any point $z \in 
X^{(n)}$ 
is of this form for a unique partition
 of $n$.  If necessary, we denote 
the partition associated with $z$  by $\bnu (z)$.

\noindent 
Let $_{x_j}\D$ be a collection of open neighborhoods of the points $x_j \in X$
subject to the following two conditions:

\smallskip
1) they are pairwise disjoint;

2) each $_{x_j}\D$ is  biholomorphic to a  unit bi-disk $\D \subseteq  \comp^2$
by a fixed isomorphism $f_{x_j}: \D \to {_{x_j}\D}$.

\smallskip
\noindent
Let $0 < \e \leq 1$ be a real number. Consider  $\D (\e):=
\{ \, (z_1, z_2) \in \D  \; : \; |z_i|  < \e , \, i=1,\, 2 \, \}$ and define
${_{x_j}\D} (\e):= f(x_j)(\D(\e ))$: it is a  neighborhood of $x_j$ in $X$.

\noindent
By  abuse of notation, we denote the $m$-fold symmetric products
$ ({_{x_j}\D} (\e))^{(m)} $ by
${_{x_j}\D}^{(m)} (\e)$.
By virtue of the defining property of symmetric products in the category
of complex spaces, the cartesian products
$$
\prod_{j=1}^k{_{x_j}\D^{(\nu_j)} (\e)}
$$
are  naturally biholomorphic to  neighborhoods of $z$ in $X^{(n)}$.
We denote them by $U_z(\e)$.

\smallskip
\noindent
If $z$ is fixed, then  the open sets $U_z(\e)$, $0< \e \leq 1$,
 form a fundamental system of
 neighborhoods of
$z$ in $X^{(n)}$. 

\noindent
If $z$ varies as well, then the open sets $U_z(\e)$ 
form a basis for the topology of $X^{(n)}$. We call these neighborhoods {\em basic}.

\smallskip
\noindent
Note that: 1) if $z \in X^{(n)}_{(\bnu)}$, then $U_z(1)$ does not meet
any stratum corresponding to a partition $\tilde\mu \neq \bnu$ such that
$\lambda (\bnu) \geq \lambda (\tilde{\mu})$, and 2) if  
$U_z(1) \cap U_w(1) \neq \emptyset$, then  either $\bnu(w) \geq \bnu (z)$ 
(i.e.
$
X^{(n)}_{( \bnu(z) )}
\subseteq \overline{X^{(n)}_{( \bnu(w) )}}
$), or $\bnu (z) \geq \bnu (w)$.

\subsection{The spaces $\XA$ and the morphisms $K_{\ba{a}}$}
\label{thespaces}
Let $\ba{a} \in P(n)$ be  a fixed partition of $n$ in $a$-notation.
Set $X^{(0)}:=pt$, where $pt$ is a single fixed point viewed
as a complex space. Define
$$
\XA:= \prod_{i=1}^n{ X^{(a_i)} }.
$$
By abuse of notation, a  point in $\XA$
will be denoted by
an $n$-tuple
$$
\tau:= (\, x_1^1+ \ldots + x^1_{a_1} \, , \,  \ldots \, ,
\,  x_1^n+  \ldots + x^n_{a_n}\, )
$$
where it is understood that for every index $i$ for which
$a_i=0$, the corresponding
entry
$(\, x_1^i+ \ldots + x^i_{a_i}\, )$ is to be replaced by  the point $pt$.

\noindent
By the defining property of symmetric products in the category
of complex spaces, there is a morphism
$$
K_{\ba{a}}: \XA \longrightarrow X^{(n)},
$$
$$
(x_1^1+ \ldots + x^1_{a_1}, \ldots, x_1^n+  \ldots + x^n_{a_n})
\longrightarrow
\sum_{i=1}^n{i(x_1^i+ \ldots + x^i_{a_i})},
$$
where, by abuse of notation,  for every
index $i$ for which $a_i=0$ the  corresponding summand in the 
sum on the right is omitted.

\medskip
The image of $K_{\ba{a}}$
is  the closure $\overline{\Xa}$ of the
stratum associated with $\ba{a}$.
\begin{lm}
\label{normal}
 The morphism
$$
K_{\ba{a}}: \XA \longrightarrow  \overline{\Xa}
$$
is the normalization map.
\end{lm}
{\em Proof.}
The morphism $K_{\ba{a}}$ is proper, finite and bimeromorphic
onto its image.
Since $\XA$ is a normal complex space, it is the normalization of
$\overline{\Xa}$.
\blacksquare

\medskip
\noindent
Define
$$
K: = \coprod_{\ba{a} \in P(n)} K_{\ba{a}} : \coprod_{\ba{a} \in P(n)} \XA
\longrightarrow X^{(n)},
$$ 
and, for every integer $1\leq l\leq n$,

$$
K_l: = \coprod_{\lambda(\ba{a})=l } K_{\ba{a}}: \coprod_{
\lambda (\ba{a} )=l } \XA \longrightarrow \bigcup_{
\lambda (\ba{a} )=l}{\overline{\Xa}}= \coprod_{
\lambda (\ba{b} )\leq l } X^{(n)}_{(\ba{b})} \subseteq X^{(n)}.
$$

\subsection{The sets $F_z$, $S(\bnu)$, $F_z^{\ba{a}}$, $S^{\ba{a}}(\bnu)$,  $F_z(l)$
and $S(l,\bnu)$ }
\label{thesets}

\noindent
Let $z\in X^{(n)}$ be a point and 
$\bnu:= \bnu(z)$ so that $z=\sum_j{\nu_j x_j}$ for pairwise distinct points
$x_j \in X$.
For every $\ba{a} \in P(n)$ and every  integer $1\leq l \leq n$
it is convenient to
define
$$
F_z:=K^{-1}(z); \quad F_z^{\ba{a}}:=K_{\ba{a}}^{-1}(z); \quad
F_z(l):= \coprod_{\lambda (\ba{a}) =n- l}{F_z^{\ba{a}}}.
$$
The reason for the notational switch $l \to n-l$  is the following
\begin{lm}
\label{supp} 
Let $h$ be a non-negative integer. Then:

$$
Supp (R^{2h} \pi_* \rat_{X^{[n]}}) = 
\coprod_{
\lambda (\ba{a} )\leq n-h } X^{(n)}_{(\ba{a})}
$$
(this set is empty for $h\geq n$).
$$
Supp (R^{2h+1} \pi_* \rat_{X^{[n]}}) = \emptyset.
$$
\end{lm}
{\em Proof.}
%Let $z=\sum_{j=1}^n{\nu_j x_j} \in X^{(n)}_{(\bnu)}$. Then
%$\pi^{-1}(z) \simeq \prod_{j=1}^n{ \D^{[\nu_j]}_{\nu_j {\bf 0} } }$.
The assertion follows from Theorem  \ref{pobd}.2, the local description
of $\pi$ given in \S\ref{localmodel} and   K\"unneth Formula.
\blacksquare

\bigskip
\noindent
Let $\beta=(\ba{b_1}, \ldots, \ba{b_k} ) \in \prod_{j=1}^k{P(\nu_j)}$
be a $k$-tuple of partitions in $a$-notation, and $\ba{b_j}(i)$ be  the 
$i$-th entry of the partition $\ba{b_j}$. 
Define a map 
$$u: \prod_{j=1}^k{P(\nu_j)} \longrightarrow P(n)
$$ 
as follows:
let $\ba{b_j}(i)=0$, if $i>\nu_j$ and 
assign to  a $\beta$ as above, the partition $u(\beta)$ of $n$ which has 
$\sum_{j=1}^k{\ba{b_j}(i)}$ as $i$-th entry
in $a$-notation. 
Define, for every $0\leq l \leq n-1$ and for every  $\ba{a} \in P(n)$
$$
S(\bnu):= \prod_{j=1}^k{P(\nu_j)}, \quad S^{\ba{a}}(\bnu):= u^{-1}(\ba{a}) ,
\quad S(l,\bnu):  = \coprod_{\lambda (\ba{a})=n-l}{S^{\ba{a}} (\bnu)} .
$$
\medskip
\noindent
The  set $F_z$ is the disjoint union of all sets $F_z^{\ba{a}}$.
In order to
determine all the sets above it is sufficient to 
determine $F_z^{\ba{a}}$, for every $\ba{a} \in P(n)$.

\smallskip
\noindent
Note that $F_z^{\ba{a}}\neq \emptyset$ iff $\ba{a}\geq \bnu$, i.e.
iff $X^{(n)}_{(\bnu)} \subseteq \overline{\Xa}$.

\noindent
In order to determine $F_z^{\ba{a}}$ we need to impose  the condition
$$
K_{\ba{a}}(\tau)=
\sum_{i=1}^n{i\,(\,x_1^i+ \ldots + x^i_{a_i}\,)} = \sum_{j=1}^k{\nu_j x_j}.
$$
It follows that 
$$
\tau=(\, \e_{11}x_1 + \ldots +  \e_{1k}x_k \, , \, 
 \ldots \, ,  \, \e_{n1}x_1 + 
\ldots \e_{nk}x_k \, ),
$$
where the $\e's$ are integers subject to the  conditions:

\smallskip
1)   $\e_{lm} \geq 0$, $\forall \,  l,m$;

2) for every $j=1, \ldots, k$, we have $\sum_{i=1}^n{i\, \e_{ij}} =\nu_j$;

\noindent
to introduce the third condition, note that by 1) and 2) above,
a solution $\e_{lm}$ gives an element $\beta=(\ba{b_1}, \ldots,
\ba{b_k}) \in S(\bnu)$:
for every $j=1, \ldots, k$ the sequence $(\e_{1j}, \ldots, \e_{nj})$
is a partition of $\nu_j$ in $a$-notation; the third condition
is that

3) $\beta \in S^{\ba{a}}(\bnu)$.

\bigskip
\noindent
There are natural bijections
$$
F_z \Longleftrightarrow S(\bnu), \quad F_z^{\ba{a}} \Longleftrightarrow
S^{\ba{a}}(\bnu), \quad F_z(l) \Longleftrightarrow \coprod_{\lambda (\ba{a})
=n-l}{S^{\ba{a}}(\bnu)} = S(l, \bnu).
$$
The elements of $F_z$ will be denoted by ${\zeta_z}_{\beta}$, the ones
of $F_z^{\ba{a}}$ by
${\zeta_z^{ \ba{a} } }_{\beta}$ 
and the ones of $F_z(l)$ 
by
$\zeta_z(l)_\beta$.

\medskip
\noindent
Fix an integer $0\leq h \leq n-1$. 
The set $S(h,\bnu)$ can also be described as follows.
Let $k:=\lambda (\bnu)$ and  $H$ be the set   of
$k$-tuples $h=(h_1, \ldots, h_k)$ of non-negative integers with the property
that $h= \sum{h_j}$. Then  $ S (h, \bnu) \subset \prod_{j=1}^k P(\nu_j)$
is the set of  $k$-tuples $\beta=\{ \ba{b_1}, \ldots,
\ba{b_k} \}$  of partitions ${\bf{\ba{b_j}}} \in P(\nu_j)$
 with $\lambda (\ba{b_j})=\nu_j-h_j$ for some $(h_1, \ldots, h_k) \in H$.

\subsection{Local model for the Douady-Barlet morphism}
\label{localmodel}
Let $z$ be a point in $X^{(n)}$. This point determines the unique
stratum $X^{(n)}_{(\bnu(z))}$ on which it lies.  Let
$\nu_1, \ldots, \nu_k$ be the entries of $\bnu(z)$.
Let 
$U_z(\e)= \prod_{j=1}^k{_{x_j}\D^{(\nu_j)}}(\e) $ be a basic
neighborhood of 
$z$ in $X^{(n)}$.
By abuse of notation,  we shall write
$_{x_j}\D^{[\nu_j]} (\e)$ in place of  $(_{x_j}\D (\e))^{[\nu_j]}.$

\noindent
Around $z$, the Douady-Barlet morphism 
can  and will be identified with the  natural morphism 
$$\prod_{j=1}^k{_{x_j}\D}^{[\nu_j]} (\e)
 \longrightarrow
\prod_{j=1}^k{{_{x_j}\D}^{(\nu_j)} (\e) }.
$$
In fact we can identify naturally the space on the left with
the pre-image under $\pi$ of the space on the right.

\subsection{The morphism of complexes}
\label{main construction}
We first make the following elementary observation which simplifies the
picture.
\begin{lm}
\label{nomorecomp}
 $\overline{ X^{[n]}_{(\bnu)} }=\pi^{-1}(\overline{ X^{(n)}_{(\bnu)} })$

\end{lm}
{\em Proof.}
By the local description of the Douady-Barlet map it is enough to prove
the following:

\smallskip
\n
{\em Claim:
Let ${\cal I}$ be the ideal of a subscheme of $\D$ of length $m$
concentrated at $o \in \D$. Let $\tilde{\mu}=(\mu_1, \ldots,$ $ \mu_l)$ be a partition of $m$. There exists 
a family of subspaces in the stratum $ X^{[m]}_{(\tilde{\mu})}$ specializing to
${\cal I}$}.

\n
By virtue of Brian\c{c}on's result Corollary \ref{brianc},
 $\D^{[m]}_o$ is irreducible of 
dimension $m-1$, so that it is enough to prove the {\em Claim} for every ${\cal I}$
in a dense subset $U$ of $\D^{[m]}_o$.
We choose $U$ to be the subset of curvilinear subspaces in
$\D^{[m]}_o$. After a change of variables, such an ideal can be written as
${\cal I}= (y- c_1 x - \ldots - c_{n-1}x^{n-1}, \, x^n)$, $c_i \in \comp$. The seeked for family is given by
$(y- c_1x - \ldots -c_{n-1}x^{n-1}, (x - \gamma_1)^{\mu_1}
\ldots (x-\gamma_l)^{\mu_l})$, $\gamma_i \in \comp$.
\blacksquare

\medskip
As a consequence, $\pi_{| \overline {X^{[n]}_{(\bnu)}}}:
\overline{ X^{[n]}_{(\bnu)} }=\pi^{-1}(\overline{ X^{(n)}_{(\bnu)} })
 \longrightarrow \overline{X^{(n)}_{(\bnu)}}$
 has irreducible fibers and, for any $U \subseteq X^{(n)}$ and any partition $\bnu$,
 the irreducible components of  $\overline {X^{[n]}_{(\bnu)}} \cap \pi^{-1}(U)$ are in one to one correspondence with those of
$\overline { X^{(n)}_{(\bnu)}} \cap U$. Note that if we did not prove
Lemma \ref{nomorecomp} we would still have a natural injective correspondence
between irreducible components on the target and on the source of $\pi$.

\bigskip
To perform our basic construction it is natural to use Borel-Moore homology. 
We recall  a few basic facts about it ( cfr. \ci{ful}, 19.1 and the references there).

a) for an n-dimensional complex space $X$, a basis for $H^{lf}_{2n}(X,\rat)$ is given by the irreducible components of $X$.

b) $H^{lf}_{\bullet}$ is covariant with respect to proper morphisms,
e.g. closed immersions.   

c) an open imbedding $j:U \longrightarrow X$ gives a restriction morphism 

$$
j^*: H^{lf}_{\bullet}(X,\rat) \longrightarrow H^{lf}_{\bullet}(U,\rat ).
$$
Fix a partition $\bnu$ of $n$ and consider
the sheaf  $ {\cal F}^{[n]}_{\bnu}$ on $X^{(n)}$ associated with the presheaf
$$
\begin{array}{lll}
U &\longrightarrow & \{
\rat -\hbox{vector space generated by the irreducible components of } \pi^{-1}(U)\cap  \overline{X^{[n]}_{(\bnu)}} \}  \\
& = & H^{lf}_{2(n+\lambda(\bnu))}( \overline{X^{[n]}_{(\bnu)}}\cap \pi^{-1}(U), \rat)
\end{array}
$$
where the presheaf structure is given by the restriction morphisms
stemming from  c).

The push-forward associated with the closed embedding 
$ \overline{X^{[n]}_{(\bnu)}} \longrightarrow X^{[n]}$ and  Poincar\'e Duality  on $X^{[n]}$ identify   canonically 
 $ {\cal F}^{[n]}_{\bnu} $ 
with a subsheaf $R^n_{\bnu}$ of $R^{ 2(n-\lambda(\bnu)) } \pi_* \rat_{X^{[n]}}$.

From the previous discussion it follows that the presheaf is isomorphic 
to the analogous one defined by the irreducible 
components of $U\cap \overline{X^{(n)}_{(\bnu)}}$, which, by Zariski Main Theorem, is isomorphic to 
${K_{\ba{a}}}_* \rat_{X^{(\ba{a})}}$, where $\ba{a}$ is  $\bnu$ in  a-notation.

Let      $({\cal D}^{\bullet}, d)$ be the resolution of the constant sheaf 
$\comp_{X^{[n]}}$ given by the complex of currents. It is an acyclic
resolution with respect to $\pi_*$. 
Let ${\cal Z}^{\bullet}$ denote the subcomplex of closed currents. 
Every analytic cycle defines the closed current of integration along itself.
We can therefore define a morphism of sheaves
$ \Psi^{[n]}_{\bnu}: {\cal F}^{[n]}_{\bnu}\longrightarrow 
\pi_*{\cal Z}^{2(n-\lambda(\bnu))} \subseteq  \pi_*{\cal D}^{2(n-\lambda(\bnu))}$ 
whose projection on the cohomology sheaves
${\cal H}^{2(n-l(\nu))}(\pi_* {\cal D}^{\bullet},d)=
R^{2(n-\lambda(\bnu))}\pi_*\comp_{X^{(n)}} $ gives
 the previously defined identification with $R^n_{\bnu}$.

\medskip
We summarize the previous discussion in the following 
\begin {pr}
\label{prpsi}
There is a natural injective morphism of complexes of sheaves:
$$
\Psi: \,  \bigoplus_{h=0}^{n-1} \;  \bigoplus_{\lambda (\ba{a)}=n-h}
{K_{\ba{a}}}_*\comp_{X^{(\ba{a})}}[-2h] \longrightarrow \pi_* {\cal D }^{\bullet}
$$
where the l.h.s complex is endowed with zero differentials.
\end{pr}

Our main result is that $\Psi$ is a quasi-isomorphism. We shall need the following
\begin{pr}
\label{basis}
Let $U$ be a basic neighborhood of a point $z \in  X^{(n)}_{(\bnu)}$. A basis for $H^{2h}(\pi^{-1}(U), \rat)$ is given
by the cohomology classes of the irreducible components of $ X^{[n]}_{(\tilde{\mu})}\cap \pi^{-1}(U)$, for all $\tilde{\mu} \leq \bnu$,
(i.e.  $X^{[n]}_{(\bnu)} \subseteq \overline {X^{[n]}_{(\tilde{\mu})}}$) with $\lambda(\tilde{\mu})=n-h$.
\end{pr}
{\em Proof.} Since $U$ is basic, we have
$\pi^{-1}(U)\simeq \prod_{j=1}^k{_{x_j}\D}^{[\nu_j]} (\e)$
and
$$
H^{2h}(\prod_{j=1}^k{ {_{x_j}\D}^{[\nu_j]} (\e),\rat)} = 
\bigoplus_{ \sum_j{h_j}=h;\;h_j\geq 0 } \;\;
\bigotimes_{j=1}^k{ H^{2h_j} ({_{x_j}\D}^{[\nu_j]} (\e),\rat)}. 
$$
By virtue of this  K\"unneth decomposition and of Theorem \ref{pobd}.3, we can form a basis for the 
vector space above by taking  products of the subvarieties
that  we obtain via Nakajima's construction on each factor
$_{x_j}\D^{[\nu_j]}$. More precisely,
we take the cohomology classes associated with the closed subvarieties
$$
_z\D^{[\bnu]}_{\beta} (\e):= 
\{ \prod_{j=1}^k
\overline {{ _{x_j}{\D}^{[\nu_j]}_{(\bf{\ba{b_j}})}(\e) }} \}
$$
which are indexed by the $k$-tuples of partitions
$\beta=(\ba{b_1}, \ldots, \ba{b_k} ) \in  S(h, \bnu) \subseteq \prod_{j=1}^k{P(\nu_j)}$ since
$$
_z\D^{[\bnu]}_{\beta} (\e)= \overline{\pi^{-1}(
\{ \prod_{j=1}^k{ _{x_j}{\D}^{(\nu_j)}_{(\bf{\ba{b_j}})}(\e) } \})}.
$$
By virtue of the combinatorics previously developed in Section
\ref{thesets} and above, and because of Lemma  \ref{nomorecomp}, these are precisely 
the branches of  $\overline {X^{[n]}_{(\tilde{\mu})}\cap \pi^{-1}U}$ 
for $\tilde{\mu} \leq \bnu$ and $\lambda (\tilde{\mu}) =n-h$.
\blacksquare

\bigskip
\noindent
Note that the complex spaces $\sbcpt$ are irreducible in $U_z(\e)$, but not locally so: they may become locally reducible
around  certain points  $y \in U_z(\e)$; in this case $\sbcpt$
breaks up around $y$ into the union of its irreducible components.
\begin{ex}
{\rm 
 Let $n=4$, $z=4z_1 \in X^{(4)}_{(4)}$, $y=2y_1 + 2y_2
\in X^{(4)}_{(2,2)}$, where $y_1$ and $y_2$ are ``near"
$z_1$. The space  $\overline{ X^{(4)}_{(2,1,1)} }$ is locally irreducible
around $z$: there is the only branch 
$\overline{_{z_1}{\D^{(4)}_{(2,1,1)}}(\e)}$
corresponding to $\beta=(1^2,2^1)$. The closed subvariety $\overline{ X^{(4)}_{(2,1,1)} }$ is locally reducible around $y$:
there are two branches $\overline{{_{y_1}{\D^{(2)}_{(2)}}(\e)}} \times   
\overline{{_{y_2}{\D^{(2)}_{(1,1)}}(\e)}}$ and 
$\overline{{_{y_1}{\D^{(2)}_{(1,1)}}(\e)}} \times   
\overline{{_{y_2}{\D^{(2)}_{(2)}}(\e)}}$ and they correspond to 
$\beta_I=(2^1, 1^2)$ and $\beta_{II}=(1^2,2^1)$, respectively.
}
\end{ex}

\section{Decomposition Theorem
for the Douady-Barlet morphism}
\subsection{Proof of the Decomposition Theorem}
\label{decompo}
\begin{tm}
\label{decel}
Let $X$ be a complex surface.
The injective morphism of complexes  of Proposition \ref{prpsi}
$$
\Psi: \,  \bigoplus_{h=0}^{n-1} \;  \bigoplus_{\lambda (\ba{a)}=n-h}
{ K_{\ba{a}} }_* \comp_{X^{(\ba{a})}}[-2h] \longrightarrow \pi_* {\cal D}^{\bullet}
$$
is a quasi-isomorphism, i.e. it induces isomorphisms on the cohomology sheaves.
\noindent
In particular, $\real \,  \pi_* \comp_{X^{[n]}}$  is isomorphic, 
in the derived category,
to a complex with trivial differentials. 
\end{tm}
{\em Proof.}
Since the differentials of the complex on the left hand side, abbreviated
by  ``$CL$", are trivial,
 we have
\[ {\cal H}^{t}(CL) =
\left\{
\begin{array}{ll}
\bigoplus_{\lambda (\ba{a)}=n-h}
{ K_{ \ba{a} } }_*\comp_{X^{(\ba{a})}}, & 
\mbox{ if $t=  2h$ and $0\leq h \leq n-1$,}  \\
0, &  \mbox{ if  otherwise.}
\end{array}
\right. \]
The cohomology sheaves of the complex on the right hand side $CR$ are
 \[ {\cal H}^{t}(CR) =
\left\{
\begin{array}{ll}
Ker (\pi_*d^{2h})/Im (\pi_*d^{2h-1}) \simeq R^{2h}\pi_*\comp_{X^{[n]}}, 
& \mbox{ if $t=  2h$ and $0\leq h \leq n-1$,}  \\
0, &  \mbox{ if  otherwise.}
\end{array}
\right. \]
The conclusion of the theorem
is equivalent to showing that the induced map on the stalks
$({\cal H}^{ 2h}(\Psi))_z$ is bijective for every $z \in X^{(n)}$ and for
every $0\leq h \leq n-1$. This is precisely the content of  Proposition \ref{basis}
\blacksquare

\subsection{Some remarks}
\label{remarks}
\begin{rmk}
\label{meaning}
{\rm The complex $\real \, \pi_* \comp_{X^{[n]}}$ is defined up to isomorphism in 
the derived category of the category of complexes of sheaves on $X^{(n)}$.
Note that the trivial complexes ${K_{\ba{a}}}_*
\rat_{X^{(\ba{a})}}$ are a natural realization of $IC^{\bullet}(\rat_{  X^{(n)}_{(\ba{a})} })$.
The morphism $\Psi$ of 
Theorem \ref{decel} is an explicit, injective and natural quasi-isomorphism of complexes of sheaves
between a  complex with trivial differentials and a natural 
realization
of $\real \, \pi_* \comp_{X^{[n]}}$. 
This is why we call
Theorem \ref{decel} ``The Decomposition Theorem."
}
\end{rmk}

\begin{rmk}
\label{pervic}
{\rm
The authors of \ci{go-so} place themselves in the algebraic context
in order to use the Decomposition Theorem
in \ci{b-b-d}. The restriction to the algebraic case is not necessary.
Recall that the Douady-Barlet morphism is projective
by virtue of Theorem \ref{dcm}.
 M. Saito
\ci{saito}, using his theory of mixed Hodge modules,
has proved the necessary result in the analytic category
for projective morphisms. This approach
gives a not necessarily natural isomorphism in a derived category, whereas 
Theorem \ref{decel} gives an explicit quasi-isomorphism of complexes.
Our approach by-passes the use of the deep decomposition theorems
\ci{b-b-d} and \ci{saito}.
}
\end{rmk}

\begin{rmk}
{\rm
In order to use Saito's Decomposition Theorem, one  needs the irreducibility result of Brian\c{c}on \ci{br} V. 3.3  to identify the intersection
cohomology complexes  occuring in Saito's Decomposition Theorem;
see \ci{naknotes}, \S6.1 and \S6.2.
Our proof of Theorem \ref{decel} does not depend in any essential way
on Brian\c{c}on's result;  see  Corollary \ref{brianc}.
}
\end{rmk}

\begin{rmk}
 {\rm  Because of the previous remarks, our proof of Theorem \ref{decel} is different in spirit and in detail
from the one outlined above. 

\n
In addition, our   proof
of G\"ottsche Formula given below, being based on Theorem \ref{decel}, besides working also in the
analytic context, is significanlty different
from the ones in the literature concerning algebraic surfaces.
In particular, we have not used \ci{e-s} Theorem  1.1.(iv), as in \ci{go}
or \ci{ch}.
}
\end{rmk}

\begin{rmk}
\label{rational}
{\rm 
By taking a suitable resolution of $\rat_{X^{[n]}}$, one can prove, in the same way,
that there is an analogous  {\em natural} decomposition for 
$\real \, \pi_*  \rat_{ X^{[n]} }$. We do not need this fact and, for now, we omit the details.
}
\end{rmk}

\section{Consequences of the Decomposition Theorem}
\label{conseq}
\subsection{The Leray Spectral Sequence for the pair $(\pi, \rat_{X^{[n]}})$}
\label{lls}
An immediate consequence of  Theorem \ref{decel},
 is the following 
\begin{cor}
\label{isr} For every integer 
$h$ there are natural isomorphisms
$$
R^{2h}
 \pi_* \rat_{X^{[n]}} \simeq   {K_{n-h}}_* 
\rat_{
       \coprod_{ \lambda (\ba{a}) =n-h }{ X^{ (\ba{a}) } } 
}\simeq 
\bigoplus_{ \lambda (\ba{a}) =n-h } 
 {K_{ \ba{a} }}_* \rat_{ X^{ (\ba{a}) } },
$$
and $R^{2h+1} \pi_* \rat_{X^{[n]}} =0$.
\end{cor}
{\em Proof.} For $\comp$-coefficients it  follows from  Theorem \ref{decel} by taking cohomology sheaves.

\n
For every $\ba{a}=\bnu \in P(n)$, we have a natural identification of 
${K_{\ba{a}}}_* \rat_{X^{(\ba{a})}}$  with $R^n_{\bnu}$ (see \S \ref{main construction}).  The assertion for  $\rat$-coefficients follows in view
of Proposition \ref{basis}.
\blacksquare

\begin{tm}
\label{leray}
The Leray spectral sequence for the pair
$(\pi, \rat_{X^{[n]}})$ is $E_2$-degenerate.
\end{tm}
{\em Proof.} In fact Theorem
\ref{decel} and the finiteness of the morphisms
$K_{\ba{a}}$ imply the following {\em natural}
decomposition in the derived category: 
$$
\real \, \pi_* \comp_{X^{[n]}} \simeq  \bigoplus_{j\geq 0} R^j\pi_* \comp_{X^{[n]}} [-j].
$$
The $E_2$-degeneration for $\comp$-coefficients implies the one for $\rat$-coefficients.
\blacksquare

\begin{rmk}
\label{ratleray}
{\rm 
In view of Remark \ref{pervic}, one sees that a decomposition
as in  Theorem \ref{leray}   holds
with $\rat$-coefficients, however one may prefer a
natural one. One gets one by virtue of  Remark \ref{rational}.}
\end{rmk}

\subsection{G\"ottsche's Formula and Nakajima's Interpretation}
\label{goetnak}
\begin{tm}
\label{go}
{\bf (G\"ottsche's Formula )}
Let $X$ be a complex surface.

\noindent
Then for every  $l \in \zed$ and every integer $n \in \nat$, we have a natural isomorphism:
$$
H^{l} (X^{[n]}, \rat) \simeq \bigoplus_{\ba{a} \in P(n) }{
H^{l -2n  + 2\lambda (\ba{a})} (X^{(\ba{a})}  , \rat) }.
$$
Assume in addition that $X$ has  finite Betti numbers
$b_i(X)$, $i=0, \ldots, 4$. 
Then the Douady spaces $X^{[n]}$ have  finite Betti numbers and the generating function for the Poincar\'e polynomials
 is:
\begin{equation}
\label{gofock}
\sum_{n=0}^{\infty}{
P_t(X^{[n]})\, q^n} =
\prod_{m=1}^{\infty}{
\frac{ (1+  t^{2m-1}\,q^m)^{b_1(X)} \; (1+  t^{2m+1}\, q^m)^{b_3(X)}) }{
(1-  t^{2m-2}\, q^m)^{b_0(X)}\;  (1-  t^{2m}\, q^m)^{b_2(X)}
\; (1-   t^{2m+2}\, q^m)^{b_4(X)}} 
}.
\end{equation}
\end{tm}
{\em Proof.} The quasi-isomorphism of Theorem \ref{decel}.
 induces an isomorphism in   hypercohomology. 
This proves the first assertion for $\comp$-coefficients. 
The assertion
for $\rat$-coefficients follows 
%since Nakajima's 
%operators are defined over $\rat$ 
from Corollary \ref{isr} and Theorem \ref{leray}.
(cfr. also Remark \ref{rational}).
Formula (\ref{gofock}) follows from a formal manipulation which builds on Macdonald's
Formula; see \ci{naknotes}, page 69.
\blacksquare

\begin{rmk}
\label{natt}
{\rm Theorem  \ref{goetnak} is new for $X$ non algebraic. If $X$ is algebraic,
then this result is slightly  more precise than the corresponding statement
in \ci{go-so} since it distinguishes a natural
isomorphism.}
\end{rmk}
\bigskip
The following corollary has been proved independently by the first author
using elementary topology in \ci{de}. 
The r.h.s. exhibits modular behavior
\ci{go}.
\begin{cor}
\label{mark}
(Cfr. \ci{de}) Let $X$ be a complex  surface with 
finite Betti numbers. Let $e(X^{[n]})$ be the Euler number of $X^{[n]}$.
We have  the following
generating function for $e\left(X^{[n]}\right)$:
\begin{equation}
\label{en}
\sum_{n=0}^{\infty}{e\left(X^{[n]}\right)\,q^n} = \prod_{m=1}^{\infty} \left( \frac{1}{1-q^m}
\right)^{e(X)}.
\end{equation}
\end{cor}
{\em Proof.} Set $t=-1$ in G\"ottsche's formula Theorem \ref{go}.
\blacksquare
\begin{tm}
\label{naka}
({\bf Nakajima's Interpretation of G\"ottsche's Formula.})
Let $X$ be a complex surface with finite Betti numbers.
The super 
vector space
${\Bbb H}(X)= \oplus_{n\geq 0} H^*(X^{[n]},\rat)$
 is an irreducible highest weight representation as the representation of 
the Heisenberg/Clifford algebra, with highest weight
vector the generator   ${\bf 1} \in H^0(X^{[0]},\rat)\simeq \rat$. 
\end{tm}
{\em Proof.}  
By virtue of Theorem
\ref{nakgo}, Theorem \ref{go} was the only missing piece.
\blacksquare

\medskip
\begin{rmk}
\label{lehnaz}
{\rm 
For a  geometric  action of the Virasoro algebra  see \ci{le}.
}
\end{rmk}

\subsection{The Hodge structure in the K\"ahler case}
For the following theorem recall that the cohomology of the quotient 
by  a finite group of a smooth compact K\"ahler manifold carries a pure
Hodge structure.

\begin{tm}
\label{mhsxn}
Let $Y$ be a compact, K\"ahler, complex surface. Let $X$ be a Zariski-dense
open subset of $Y$. Then for every 
$l \in \zed$ and every $n \in \nat^+$  we have a natural isomorphism of mixed Hodge
structures:
$$
H^{l}(X^{[n]}) \otimes \rat (n) \simeq    \bigoplus_{\ba{a} \in P(n)}{
H^{l -2n + 2\lambda (\ba{a}) }   (X^{(\ba{a})}, \rat) \otimes \rat (\lambda (\ba{a}))}.
$$
\end {tm}
{\em Proof.} Note that $X$ has finite Betti numbers. Every class $\alpha \in H^*(\H{n})$ can be represented as
$P_{\alpha_1}[i_1]\circ \ldots P_{\alpha_r}[i_r] ({\bf 1})$.
If the classes $\alpha_i$ have  type $(p_i,q_i)$, then   $\alpha$ has type
$(\sum p_i+n-l(\tilde \nu),\sum q_i+n-l(\tilde \nu))$.
The statement now follows after a  formal manipulation.
\blacksquare 

\begin{rmk}
\label{credit}
{\rm 
In the quasi-projective case Theorem \ref{mhsxn}
(without the naturality assertion) was first proved in \ci{go-so} using Saito's Theory
of mixed Hodge modules. The remark that a proof in the projective case
depending on Nakajima's operators 
is possible can be found in  \ci{gr} and has also been made by Nakajima in a private communications to us.

\n
Note that one can compute the generating function for the virtual Hodge numbers
as in \ci{ch}.}
\end{rmk}

\subsection{Connection with Equivariant $K$-theory}
\label{eqk}
In this section we prove Theorem \ref{eqco}.

Our original motivation was to explain the following
three sets of equalities by means of the existence of {\em natural} isomorphisms.

\n
Let $G\subseteq U(2)$ be a finite subgroup, $\comp^2/G$ be the associated ``Kleinian singularity"
and $Y$ be its associated minimal resolution of singularities.
\begin{itemize}

\item
The papers \ci{a-s} and \ci{h-h} remark that the ``orbifold Euler number"
$e(\comp^2, G)$ and the Euler number $e(Y)$ coincide.

\item
The paper \ci{h-h} remarks that the equality $e(X^n,\SS_n)=
e(X^{[n]})$ holds for any
  smooth algebraic
surface $X$.

\item
The paper \ci{wa} remarks that $\dim_{\comp}{K_{\SS_n}(Y^n) \otimes_{\zed} \comp }
=$ $\dim_{\comp}{ K(Y^{[n]}) \otimes_{\zed} \comp}$. The same is true for 
every surface $Y$.

\end{itemize}
We explain these equalities via the natural isomorphism in Theorem \ref{eqco}.

\bigskip
Let us recall some  notions and facts relating them.

\bigskip
Let $G$ be a finite group and $Y$ be a locally compact, Hausdorff and paracompact
left $G$-space.

Denote by $K_G(Y)$ the Equivariant $K$-Theory of the pair $(Y,G)$;
see  \ci{a-s}, \ci{b-c} and \ci{wa}. It is a $\zed_2$-graded abelian group.
Its ``even" part is the Grothendieck group generated by $G$-vector bundles.
We shall not consider the multiplicative structure induced by the
tensor product.

Let $a\in G$ be an element and define $Y^a:= \{y \in Y \, | \; ya=y     \}$.

Define $\hat{Y} := \{ (y,b)\in Y\times G, \quad yb=y \}$.
There is a natural identification $\hat{Y}= \coprod_{g\in G} Y^g$.
There is a natural $G$-action of $G$ on $\hat{Y}$ given by
$(y,c)d:= (yd,d^{-1}cd)$.

Let $G_*$ be the set of conjugacy classes of $G$, $g \in G$ be an element, 
$[g]$ be its conjugacy class, 
 and $Z_G(g)$ be the centralizer
of $g$ in $G$; this subgroup acts on $Y^g$. If $[g]=[g']$, then there is a canonical identification
$Y^g/Z_G(g) \simeq Y^{g'}/Z_G (g' )$.

 Choose representatives
$\ba{g}=\{g_1, \ldots , g_{|G_*|} \}$ for each conjugacy class in $G_*$.

There is a homeomorphism  $\alpha_{\ba{g}}: \hat{Y}/G \simeq \coprod_{[g] \in G_*}{Y^{g_l}/Z_G(g_l)}$.

The following relates the $G$-Equivariant $K$-Theory of $Y$ to the $K$-Theory
of the fixed-point-sets.

\begin{tm}
\label{as}
(Cfr. \ci{a-s}, \ci{b-c}.)
Let $G$, $Y$, $\hat{Y}$, $\ba{g}$ and $\alpha_{\ba{g}}$ be as above. There are $\zed_2$-graded  isomorphisms
of $\comp$-vector spaces:
$$
\phi_{\ba{g}}: K_G(Y)\otimes_{\zed} \comp \longrightarrow
K(\hat{Y}/G)\otimes_{\zed} \comp \stackrel{ \alpha_{\ba{g}}^{-1} }\longrightarrow
\bigoplus_{l=1}^{|G_*|}{K ( Y^{g_l}/Z_G(g_l) ) \otimes_{\zed} \comp},
$$
where the first one is natural and the second one depends on $\ba{g}$.
\end{tm}

\begin{rmk}
\label{qnotc}
{\rm 
If $G=\SS_n$, then Theorem \ref{as} holds with $\rat$-coefficients.
}
\end{rmk}

\begin{tm}
\label{eqco}
{\bf (Connection with Equivariant $K$-theory)}
Let $X$ be a complex surface.
For every natural integer $n$ there is a natural $\zed_2$-graded 
$\rat$-linear 
isomorphisms:
$$
K_{\SS_n}(X^n) \otimes_{\zed} \rat  \simeq
K(X^{[n]}) \otimes_{\zed} \rat.
$$

\end{tm}
{\em Proof.}
Let $Y:=X^n$ and $G:= \SS_n$, where the action is given by the permutation of the factors.
There is a natural and well-known identification $P(n)={\SS_n}_*$. Having made a choice $\ba{g}$
as above we obtain, as in \ci{h-h}, a natural identification
$\beta: \hat{Y}/\SS_n \simeq  \coprod_{\ba{a} \in P(n)}{X^{(\ba{a})}}$ which does not depend on $\ba{g}$.

\noindent
By virtue of Theorem \ref{as}, of the existence of $\beta$ and of Theorem \ref{go}, we get
a natural $(\zed /2\zed)$-graded isomorphism of graded $\rat$-vector spaces:
$$
K_{\SS_n}(X^n)\otimes_{\zed} \rat \simeq K(\hat{Y}/\SS_n)\otimes_{\zed} \rat
\simeq \bigoplus_{   \ba{a} \in P(n)  }{  K( X^{ (\ba{a}) } )       \otimes_{\zed} \rat}
\simeq K(X^{[n]}) \otimes_{\zed} \rat,
$$
where the last isomorphism stems from Theorem 
\ref{go} after taking the Chern Character isomorphism:
$K(-)\otimes_{\zed} \rat \simeq  H^*(-, \rat)$.
\blacksquare
\begin{rmk}
{\rm
\bigskip
After we proved Theorem \ref{eqco}, we  received a copy
of \ci{b-g} which contains
a    statement similar to the one of  Theorem  \ref{eqco}, but where
the natural map proposed in \ci{b-g} should be constructed in an entirely different way. 
We thank V. Ginzburg for giving us a copy of \ci{b-g}.
}
\end{rmk}
\begin{rmk}
\rm{
In a lecture at Cambridge  G. Segal \ci{se} (see also \ci{wa}) introduced new structures on the 
vector space ${\Bbb K}(X):= \oplus_{n\geq 0}{K_{\SS_n}(X^n)}\otimes_\zed \rat$. }
\end{rmk}

%\noindent
%By Atiyah-Segal 
%(the symmetric groups is special: the pairing irreducible 
%representation-conjugacy classes give an invertible matrix with integral entries; the character is obtained by inverting this matrix and therefore is rational; even for the cyclic groups the characters will be 
%roots of unity; in general should be in $ \overline{\rat}$)
%we have a natural decomposition induced by the {\em character}
%$$
%K_{\SS_n}(X^n)\otimes \rat  \simeq \bigoplus_{[g]}{ K( (X^n)^g )^{C_g} \otimes \rat}
%$$
%where $[g]$ ranges over the set of   conjugacy classes
%of $\SS_n$ and $C_g$ is the centralizer of $g$.

%\bigskip
%\noindent
% $1991$ {\em Mathematics Subject Classification}: . 

%\smallskip
%\noindent
%{\em Keywords and phrases}:  .

\bigskip
\noindent
Author's address:

\smallskip
Mark Andrea A. de Cataldo,
Department of Mathematics,
Harvard University,
One Oxford Street, Cambridge, MA 02138, USA. \quad 
e-mail: {\em mde@abel.math.harvard.edu}

\smallskip
Luca Migliorini,
Dipartimento di Matematica Applicata "G.Sansone",
Via S.Marta, 3,
50139 Firenze, ITALY. \quad
e-mail: {\em luca@poincare.dma.unifi.it}

\end{document}